%% file: MicGen2.tex
\documentclass[12pt]{amsart}
\usepackage{amssymb,amsthm}
\usepackage{amsmath}
\usepackage{tikz}
\usepackage{graphicx}
\usepackage{subcaption}
\usepackage[shortlabels]{enumitem}
\usetikzlibrary{arrows}
\usepackage{float}
\usepackage{graphicx}
\usepackage{subcaption}
\usepackage{hyperref}
\usepackage{pgfplots}

\providecommand{\keywords}[1]
{
  \small	
  \textbf{\textit{Keywords---}} #1
}

\graphicspath{{AnalBragg/}}
\usetikzlibrary{decorations.pathmorphing}
\tikzset{snake it/.style={decorate, decoration=snake}}

\makeatletter
\newcommand*{\rom}[1]{\expandafter\@slowromancap\romannumeral #1@}
\makeatother

\numberwithin{equation}{section}

\theoremstyle{plain}
\newtheorem{theorem}{Theorem}
\numberwithin{theorem}{section}

\newtheorem{corollary}[theorem]{Corollary}
\theoremstyle{definition}
\newtheorem{definition}[theorem]{Definition}
\theoremstyle{remark}
\newtheorem{remark}[theorem]{Remark}
\newtheorem{example}[theorem]{Example}
\theoremstyle{remark}

\theoremstyle{remark}

\newcommand{\smo}{\setminus \mathbf{0}}
\newcommand{\zero}{\mathbf{0}}

\newcommand{\norm}[1]{\left\lVert#1\right\rVert}      % Norm
\newcommand{\abs}[1]{\left|#1\right|}                 % Absolutbetrag
\newcommand{\paren}[1]{\left(#1\right)}               % Klammern
\newcommand{\bparen}[1]{\left[#1\right]}               % eckige Klammern
\newcommand{\sparen}[1]{\left\{#1\right\}}      % Mengenklammer
\newcommand{\dd}{\mathrm{d}}  % without the space
\newcommand{\vp}{\varphi}

\newcommand{\Ac}{\mathcal{A}}

\newcommand{\Cc}{\mathcal{C}}

\newcommand{\Dc}{\mathcal{D}}
\newcommand{\Ec}{\mathcal{E}}
\newcommand{\Fc}{\mathcal{F}}

\newcommand{\Lc}{\mathcal{L}}

\newcommand{\Sc}{\mathcal{S}}

\newcommand{\WF}{\mathrm{WF}}                         % Wavefront set
\newcommand{\wf}{\mathrm{WF}}                         % Wavefront set

\newcommand{\inv}{^{-1}}

\newcommand{\partyf}[2]{\frac{\partial #2}{\partial y_{#1}}}

\newcommand{\vu}{{\mathbf{u}}}

\newcommand{\vw}{{\mathbf{w}}}

\newcommand{\bpm}{\begin{pmatrix}}
\newcommand{\epm}{\end{pmatrix}}

\newcommand{\xo}{x_0}
\newcommand{\xn}{x_n}
\newcommand{\vx}{{\mathbf{x}}}
\newcommand{\vxo}{\mathbf{x}_0}

\newcommand{\vy}{{\mathbf{y}}}

\newcommand{\vs}{\mathbf{s}}

\newcommand{\vxi}{{\boldsymbol{\xi}}}
\newcommand{\vxio}{{\boldsymbol{\xi}_0}}

\newcommand{\xin}{\xi_n}

\newcommand{\om}{\omega}
\newcommand{\vom}{\boldsymbol{\omega}}
\newcommand{\vlam}{\boldsymbol{\lambda}}

\newcommand{\veta}{{\boldsymbol{\eta}}}

\newcommand{\vsig}{{\boldsymbol{\sigma}}} %when $\sigma$ is one dimensional, 
\newcommand{\snmt}{S^{n-2}}

\newcommand{\oinf}{(0,\infty)}

\newcommand{\rr}{{{\mathbb R}}}

\newcommand{\rtwo}{{{\mathbb R}^2}}
\newcommand{\rthree}{{{\mathbb R}^3}}
\newcommand{\rn}{{{\mathbb R}^n}}

\newcommand{\dR}{\dot{\mathbb{R}}}
\newcommand{\drn}{{\dot{{\mathbb R}^n}}}
\newcommand{\rnmo}{{\mathbb{R}^{n-1}}}

\newcommand{\st}{\hskip 0.3mm : \hskip 0.3mm}

\newcommand{\be}{\begin{equation}}
\newcommand{\ee}{\end{equation}}
\newcommand{\bea}{\begin{eqnarray}}
\newcommand{\eea}{\end{eqnarray}}
\newcommand{\bean}{\begin{eqnarray*}}
\newcommand{\eean}{\end{eqnarray*}}

\newcommand{\bel}[1]{\begin{equation}\label{#1}}

\newcommand{\eel}[1]{{\label{#1}\end{equation}}}
\newcommand{\intt}{{\operatorname{int}}}

\newcommand{\noj}{{(-1)^j}}
\newcommand{\nojp}{{(-1)^{j+1}}}
\newcommand{\nojxxo}{{(-1)^j(x_1-\xo)}}
\newcommand{\pilj}{\Pi_L^{(j)}}
\newcommand{\loc}{{\text{loc}}}
\usepackage{fullpage}

\usepackage{datetime}
\usetikzlibrary{quotes, angles}

\title{Microlocal analysis of generalized Radon transforms from
scattering tomography \\{\footnotesize\ddmmyyyydate\today~\currenttime}}
\author{James W. Webber}
\address{Department of Electrical and Computer
Engineering, Tufts University, Medford, MA USA}
\email{James.Webber@tufts.edu}
\author{Eric Todd Quinto}
\address{Department
of Mathematics, Tufts University, Medford, MA USA}
\email{Todd.Quinto@tufts.edu}
\begin{document}

\begin{abstract}
Here we present a novel microlocal analysis of generalized Radon
transforms which describe the integrals of $L^2$ functions of compact
support over surfaces of revolution of $C^{\infty}$ curves $q$. We
show that the Radon transforms are elliptic Fourier Integral Operators
(FIO) and provide an analysis of the left projections $\Pi_L$. Our
main theorem shows that $\Pi_L$ satisfies the semi-global Bolker
assumption if and only if $g=q'/q$ is an immersion. An analysis of the
visible singularities is presented, after which we derive novel
Sobolev smoothness estimates for the Radon FIO. Our theory has
specific applications of interest in Compton Scattering Tomography
(CST) and Bragg Scattering Tomography (BST). We show that the CST and
BST integration curves satisfy the Bolker assumption and provide
simulated reconstructions from CST and BST data.  Additionally we give
example ``sinusoidal" integration curves which do not satisfy
Bolker and provide simulations of the image artefacts.  The observed
artefacts in reconstruction are shown to align exactly with our
predictions.
\end{abstract}


\maketitle
%\titlep
%\tableofcontents

\input{Intro}
\input{Defns}

\input{MainThm.tex}

\input{R2_2.tex}

\input{Rn.tex}
\input{Sobolev.tex}

\input{ExamplesR2_2.tex}

\section{Conclusions and further work} Here we have presented a novel
microlocal analysis of a generalized cone Radon transform $R$, which
defines the integrals of $f\in
L^2_c(\mathbb{R}^{n-1}\times(0,\infty))$ over the $(n-1)$-dimensional
surfaces of revolution of smooth curves $q$. We proved that $R$ is an
elliptic FIO order $\frac{1-n}{2}$, and we gave an explicit expression
for the left projection $\Pi_L$. Our main theorem (Theorem
\ref{thm:BolkerRn}) shows that $\Pi_L$ satisfies the semi-global
Bolker assumption if and only if $g=q'/q$ is an immersion. Two main
applications of this theory are in Compton camera imaging in CST, and
crystalline structure imaging in BST and airport baggage screening.
In section \ref{sect:ex} we showed that the CST and BST integration
curves satisfied the conditions of Theorem \ref{thm:BolkerRn}, thus
proving that the CST and BST Radon FIO satisfy the Bolker assumption.
Additionally we gave example ``sinusoidal" $q$ in example \ref{ex4}
for which the corresponding Radon transforms violate the Bolker
assumption, and we provided simulated image reconstructions from
sinusoidal Radon data. We saw artefacts appearing along the sinusoidal
curves which intersected the singularities of $f$ normal to the
direction of the singularity. The artefacts observed in reconstruction
were shown to align exactly with our predictions and the results of
Theorem \ref{thm:BolkerRn}.

The theory presented here explains some key microlocal properties of a
range of Radon transformations in $\mathbb{R}^n$, whereby the
integrals are taken over generalized cones with vertex constrained to
the $\vx'=(x_1,\ldots,x_{n-1})$ plane. In further work we aim to generalize
the set of cone vertices to suit a wider range of imaging geometries.
For example, we could consider the vertex sets which are smooth $n-1$
manifolds in $\mathbb{R}^n$, in a similar vein to
\cite{zhang2020recovery} in $\mathbb{R}^3$. 

It is noted that quality of reconstruction (with zero noise), from CST and BST data, is low using the methods considered, and there are significant boundary artefacts in the reconstructions presented (see examples \ref{ex1} and \ref{ex2}). The reconstruction methods used here were chosen to highlight the image artefacts predicted by our theory, so this is as expected. In future work we aim to derive practical reconstruction algorithms and regularization penalties to combat the artefacts, for example using smoothing filters as in \cite{vline,FrQu2015,borg2018analyzing} to remove boundary artefacts. An algebraic approach may also prove fruitful (as is discovered in \cite{WebberQuinto2020I} for CST artefacts), as this would allow us to apply the powerful regularization methods from the discrete inverse problems literature, e.g. Total Variation (TV).

\section*{Acknowledgements} We would like to thank Professor Eric Miller for his helpful suggestions, thoughts and insight towards the article, in particular towards improving the readability of the paper and helping us communicate the main ideas to a practical audience. This material is based upon work
supported by the U.S.\ Department of Homeland Security, Science and
Technology Directorate, Office of University Programs, under Grant
Award 2013-ST-061-ED0001. The views and conclusions contained in this
document are those of the authors and should not be interpreted as
necessarily representing the official policies, either expressed or
implied, of the U.S. Department of Homeland Security.  The work of the
second author was partially supported by U.S.\ National Science
Foundation grant DMS 1712207.  Similarly, the opinions, findings, and
conclusions or recommendations expressed here do not necessarily
reflect the views of the National Science Foundation.

%The authors thank the referees for thorough reviews and thoughtful
%comments that improved the article.

\bibliographystyle{abbrv}
\bibliography{RefMicGen2}

\appendix \section{Bragg curve analysis for $x_2\in(-1,1)$}
\label{appA} Throughout this section we will use the notation of
\cite{Web}, where $x_2\in(-1,1)$ now describes the coordinates of the
scanned line profile in BST (i.e., the vertex of the V is at
the point $(x_1,x_2)$, see \cite[figure 1]{Web}); $x_1\in\mathbb{R}$
plays the same role in both articles.

In \cite{Web} the authors consider a one-dimensional set of 2-D
Radon transforms, with imaging applications in BST and spectroscopy.
In example \ref{ex2} we considered the curves of integration defined
by $q_B(x_1)=\frac{x_1}{\sqrt{x_1^2+1}}$. These curves describe a
special case of the Radon transforms of \cite{Web}, where the scanned
line profile is on the centerline of the imaging apparatus,
$x_2=0$.
Here we consider the general $x_2\in (-1,1)$ case. The full set of
integration curves in BST are described by \cite[equation (4.2)]{Web}:
\begin{equation}
\label{bcur}
q_B(x_1,x_2)=\frac{1}{\sqrt{2}}\sqrt{1+\frac{x_1^2-(1-x_2^2)}{\sqrt{x_1^2+(x_2+1)^2}\sqrt{x_1^2+(1-x_2)^2}}}.
\end{equation}
Note that $q_B(x_1,0)=\frac{x_1}{\sqrt{x_1^2+1}}$. In \cite{Web} the 1-D set of Radon transforms considered take integrals over the broken-ray curves $q_B(\cdot,x_2)$ for each $x_2\in(-1,1)$. The broken-ray integrals are described by the generalized cone transform of \eqref{R2Rad} with $q=q_B(\cdot,x_2)$. Note that $q_B(\cdot,x_2)>0$ for $x_1>0$, for every $x_2\in(-1,1)$, so \ref{hyp1} is satisfied. Here we aim to show that $h_B(\cdot,x_2)=\frac{1}{g_B(\cdot,x_2)}=\frac{q_B(\cdot,x_2)}{q'_B(\cdot,x_2)}$ is an immersion for each $x_2\in(-1,1)$, thus showing that the Bolker assumption is satisfied for every scanning profile considered. As the calculation of the second order derivatives of $q_B$ in $x_1$ is cumbersome, we verify Bolker numerically, for some chosen range of $x_1,x_2$. Note that we consider the reciprocal $h_B$ of $g_B$ here to avoid division by values close to zero as $x_1\to 0$, and since $g'_B(\cdot,x_2)=0\iff h'_B(\cdot,x_2)=0$ for $g_B(\cdot,x_2) : (0,\infty)\to(0,\infty)$. 
\begin{figure}[!h]
\centering
%\begin{subfigure}{0.49\textwidth}
%\includegraphics[width=0.9\linewidth, height=5.2cm, keepaspectratio]{gimage.png}
%\subcaption{$g_B(x_1,x_2)$.}
%\end{subfigure}
\begin{subfigure}{0.55\textwidth}
\includegraphics[width=0.9\linewidth, height=6cm,keepaspectratio]{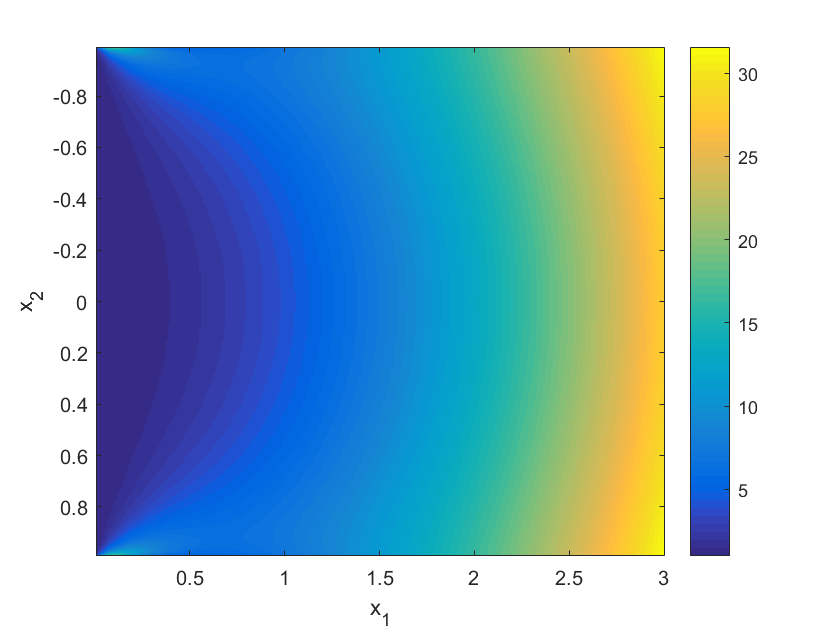}
\subcaption{$h'_B(x_1,x_2)$. $\min_{(0,3]\times(-1,1)}(h'_B)\approx 1$.}
\end{subfigure}
\caption{Numerical validation of $h'_B(\cdot,x_2)\neq 0$ for each $x_2\in(-1,1)$.}
\label{gimages}
\end{figure}
We choose to simulate $h'_B$ for $x_1\in (0,3]$ and $x_2\in(-1,1)$. The maximum $x_1$ (i.e. $x_1=3$) chosen is the maximum $x_1$ considered in the scanning setup of example \ref{ex2}. That is with scanning region $[-1,1]\times[0,2]$ and $x_0\in[-2,2]$ (the maximum $x_1$ occurs when $x_0=-2,2$ at the edge of the scanning region). We need only consider positive $x_1$ since $q_B$ of equation \eqref{bcur} is symmetric in $x_1$ about $x_1=0$. See figure \ref{gimages}, where we display $h'_B$ on $(0,3]\times(-1,1)$. Finite differences are used to approximate the derivatives. The minimum value of $h'_B$ in this range is $\min_{(0,3]\times(-1,1)}(h'_B)\approx 1$, which indicates that the Bolker assumption is satisfied in the scanning geometry of example \ref{ex2}.

\end{document}

%% file: Intro.tex
\section{Introduction}\label{sect:intro}

In this paper we present a new microlocal analysis of Radon transforms
which describe the integrals of $L^2(\mathbb{R}^n)$ functions of
compact support over the surfaces of revolution of $q\in
C^{\infty}((0,\infty))$. In 2-D the ``surface of rotation" is the
union of two curves which are the mirror images of one-another. We
denote the union of the reflected curves as ``broken-rays" (sometimes denoted by ``V-lines" in the literature \cite{vline}) when $n=2$
and the surfaces of revolution as ``generalized cones" when $n\geq 3$.
We illustrate the scanning geometry and some example integration
curves related to CST and BST in the $n=2$ case in figure \ref{fig1}. The Radon data is
$n$-dimensional and comprised of an ($n-1$)-dimensional translation by
$\vx_0\in\mathbb{R}^{n-1}$ and a one-dimensional scaling by
$E\in(0,\infty)$. We use the notation of \cite{Web} in BST, where
classically $q=\sin\theta$ denotes the sine of the Bragg angle
($\theta$) and $E$ denotes the photon energy. Here we generalize the
Radon transforms of \cite{Web} and analyze their stability
microlocally. Our theory also has applications of interest in gamma
ray source imaging in CST, specifically towards the broken-ray
transforms of \cite{ambartsoumian2012inversion,
AmbartsoumianLatifi-Vline2019,vline,florescu2010single,florescu2011inversion,
florescu2009single, morvidone2010v,truong2015new}, and the cone
Radon transforms of
\cite{
cebeiro2016back,
GouiaAmbartsoumian2014,
haltmeier2014exact,
3D2,
kuchment2016three,
maxim2009analytical,
moon2016determination,
moon2017inversion,
moon2017analytic,
nguyen2005radon,
Terzioglu-KK-cone,
3D1,
zhang2020recovery}.
\begin{figure}[!h]
\centering
\begin{tikzpicture}[scale=1]
\begin{axis}[
  legend pos=outer north east ,xlabel=$x$, ylabel=$y$, xmin=-3,
xmax=9, ymin= 0,ymax=3, xtick={-3,0,3,6,9}, ytick={0,0,1,2,3},
xticklabels={-3,0,$x_0$,6,9}, 
major tick length=0mm, grid=major
]
\addplot+[mark=none,samples=100,domain=3:9] {2*(x-3)/sqrt((x-3)^2+1)};
\addplot+[mark=none,samples=100,domain=3:9] {0.5*(x-3)};
\addplot+[mark=none,samples=100,domain=-3:3,blue] {-2*(x-3)/sqrt((x-3)^2+1)};
\addplot+[mark=none,samples=100,domain=-3:3,red] {-0.5*(x-3)};
%\addplot+[mark=none,samples=100,domain=-3:3,black,<->,thick] {0};
%\draw (1.4,0.01) circle (0pt) node [right] {$x_0$};
\legend{$y=Eq_B(x-x_0)$,$y=Eq_C(x-x_0)$}
\end{axis}
\end{tikzpicture}
\caption{The scanning geometry. The broken-ray curves displayed are
$q_C(x)=x$ and $q_B(x)=\frac{x}{\sqrt{x^2+1}}$, which are of interest
in CST and BST respectively. The curves are scaled by $E>0$,
translated by $x_0$ along the $x$ axis and reflected in the line
$x=x_0$. For example, in this case $x_0=3$, $E=2$ for $q_B$ and
$E=\frac{1}{2}$ for $q_C$.} \label{fig1}
\end{figure}

The generalized Radon transforms considered here are shown to be
elliptic FIO order $\frac{1-n}{2}$, and we give an analysis of the
left projections $\Pi_L$. Our main theorem proves that $\Pi_L$
satisfies the semi-global Bolker assumption (i.e. $\Pi_L$ is an
embedding) if and only if the quotient function $g=q'/q$ is an
immersion. Then, we consider the visible singularities in the Radon
data and provide Sobolev space estimates for the level of smoothing of
the target singularities. This serves to reduce the microlocal and
Sobolev analysis of our $n$-dimensional Radon FIO to the injectivity
analysis of the one-dimensional function $g$.

We consider two applications of our theory that are of interest,
namely Compton camera imaging in CST and crystalline structure imaging
in BST.  We show that the CST and BST integration curves satisfy the
conditions of our theorems, which, by implication, proves that
the CST and BST operators are elliptic FIO which satisfy
Bolker. Additionally we give example ``sinusoidal" $q$ for which the
corresponding transforms are shown to violate the Bolker
assumption.  In this case there are artefacts appearing along the $q$
curves at the points where $g=q'/q$ is non-injective. Using the $g$
mapping, we are able to predict precisely the locations of the
artefacts in reconstruction. To verify our theory, we present
simulated reconstructions of a delta function and a characteristic
function on a disc from CST, BST and sinusoidal data. The
predicted artifacts are shown to align exactly with those observed in
reconstruction.

The literature includes the microlocal analysis of broken-ray
transforms in \cite{vline,terzioglu2019some} and cone Radon transforms
in \cite{terzioglu2019some,zhang2020recovery}. In \cite{vline}, the
authors analyze the boundary artefacts in reconstruction from
broken-ray (denoted V-lines in \cite{vline}) integrals, which occur
along broken-ray curves at the edge of the data set. A smooth cut off
in the frequency domain is later introduced to combat the boundary
artefacts.  Proof of FIO and injectivity analysis of the $\Pi_L$ is
not considered however. We aim to cover this here for the broken-ray
transform. In \cite{zhang2020recovery}, the author considers the
five-dimensional set of cone integrals in $\mathbb{R}^3$, where the
cone vertices are constrained to smooth 2-D surfaces $\mathcal{S}$ in
$\mathbb{R}^3$. In \cite[Proposition 4]{zhang2020recovery} the normal operator of the cone transform is proven to be an elliptic Pseudo
Differential Operator (PDO) order $-2$ (under certain visibility
assumptions), thus implying that the
Bolker assumption is satisfied. In contrast, for $\rthree$, we consider the three-dimensional
subset of the Radon data when the surface of cone vertices
$\mathcal{S}=\mathbb{R}^2$ is the $(u_1,u_2)$ plane and the axis of
rotation has direction $\beta=(0,1)$ (using the notation in
\cite[Example 1]{zhang2020recovery}).  We prove that the Bolker
assumption is satisfied here with limited data, and our surfaces of
integration are more general than cones.  In
\cite[section 4]{terzioglu2019some} the $n$-dimensional case for the
cone transform is analyzed microlocally; the Radon integrals are taken
over the full set of cones in $\mathbb{R}^n$, and the data set is
$2n$-dimensional.  In \cite[Theorem 14]{terzioglu2019some} it is
proven that the normal operator of the $2n$-dimensional cone transform is a PDO. We consider
the $n$-dimensional subset of the Radon data where (using the notation
of \cite{terzioglu2019some}) $\vu\in\{u_n=0\}$ is constrained to the
$(u_1,\ldots,u_{n-1})$ plane, and the axis of rotation has direction
$\beta=(\zero,1)$. That is, we consider the vertical (i.e.
$\beta=(\zero,1)$) cones with vertices on $\{u_n=0\}$. The results of
\cite[Theorem 14]{terzioglu2019some} are not sufficient to prove
Bolker with limited data. The vertical cone Radon transform is also considered in \cite{3D1}, but no microlocal analysis is given. We aim to cover this important limited data case here.  In addition, our theorems are valid, not only for cones, but
for general surfaces of revolution satisfying \eqref{g'}. 

Our transform is a Radon transform on surfaces of revolution that are
generated by translation by directions in the $x_n=0$ plane.  Radon
transforms on surfaces of revolution have been considered in the pure
mathematical community, such as \cite{Cormack1987, Kurusa1993}, but in
those articles, the surfaces are generated by rotation about the
origin not translation in a hyperplane.

In \cite{Web} Radon models (denoted by the ``Bragg transform") are
introduced for crystalline structure imaging in BST and airport
baggage screening. The curves of integration in BST are illustrated by
$q_B$ in figure \ref{fig1}. Injectivity proofs and explicit inversion
formulae are provided for the Bragg transform in \cite[Theorem
4.1]{Web}. The stability analysis is not covered however. We aim to
address the stability aspects of the Bragg transform here from a
microlocal perspective. 

The remainder of this paper is organized as follows. In section
\ref{sect:defns} we recall some notation and definitions from
microlocal analysis which will be used in our theorems. In section
\ref{sect:main} we define the generalized cone Radon transform $R$,
which describes the integrals of $L^2$ functions over the surfaces of
revolution of smooth $q$. We prove that $R$ is an elliptic FIO order
$\frac{1-n}{2}$ and provide expression for the left projection
$\Pi_L$. We then go on to prove our main theorem, which shows that
$\Pi_L$ is an injective immersion if and only if $g=q'/q$ is an
immersion. The smoothing in Sobolev norms is later explained in
section \ref{sect:Sobolev}. In section \ref{sect:ex} we show that the
curves of integration in CST and BST (as displayed in figure
\ref{fig1}) satisfy the conditions of our theorems, and we provide
simulated reconstructions from CST and BST data. Additionally, in
example \ref{ex4}, we give example ``sinusoidal" $q$ with $g$ not
an immersion, thus violating Bolker. We simulate the artefacts in
reconstruction from these sinusoidal integrals. The observed artefacts are
shown to align exactly with our predictions and the theory of section
\ref{sect:main}.

%% file: Defns.tex
	\section{Microlocal definitions}\label{sect:defns} 

We next provide some
notation and definitions.  Let $X$ and $Y$ be open subsets of
$\rn$.  Let $\Dc(X)$ be the space of smooth functions compactly
supported on $X$ with the standard topology and let $\mathcal{D}'(X)$
denote its dual space, the vector space of distributions on $X$.  Let
$\Ec(X)$ be the space of all smooth functions on $X$ with the standard
topology and let $\mathcal{E}'(X)$ denote its dual space, the vector
space of distributions with compact support contained in $X$. Finally,
let $\Sc(\rn)$ be the space of Schwartz functions, that are rapidly
decreasing at $\infty$ along with all derivatives. See \cite{Rudin:FA}
for more information. 

\begin{definition}[{\cite[Definition 7.1.1]{hormanderI}}]
For a function $f$ in the Schwartz space $\Sc(\mathbb{R}^n)$, we define
the Fourier transform and its inverse 
as
\begin{equation}
%\begin{split}
\mathcal{F}f(\vxi)=\int_{\mathbb{R}^n}e^{-i\vx\cdot\vxi}f(\vx)\mathrm{d}\vx,
%\\
\qquad\mathcal{F}^{-1}f(\vx)=(2\pi)^{-n}\int_{\mathbb{R}^n}e^{i\vx\cdot\vxi}f(\vxi)\mathrm{d}\vxi.
%\end{split}
\end{equation}
\end{definition}

We use the standard multi-index notation: if
$\alpha=(\alpha_1,\alpha_2,\dots,\alpha_n)\in \sparen{0,1,2,\dots}^n$
is a multi-index and $f$ is a function on $\rn$, then
\[\partial^\alpha f=\paren{\frac{\partial}{\partial
x_1}}^{\alpha_1}\paren{\frac{\partial}{\partial
x_2}}^{\alpha_2}\cdots\paren{\frac{\partial}{\partial x_n}}^{\alpha_n}
f.\] If $f$ is a function of $(\vy,\vx,\vs)$ then $\partial^\alpha_\vy
f$ and $\partial^\alpha_\vs f$ are defined similarly.

  We identify cotangent
spaces on Euclidean spaces with the underlying Euclidean spaces, so we
identify $T^*(X)$ with $X\times \rn$.

If $\phi$ is a function of $(\vy,\vx,\vs)\in Y\times X\times \rr^N$
then we define $\dd_{\vy} \phi = \paren{\partyf{1}{\phi},
\partyf{2}{\phi}, \cdots, \partyf{n}{\phi} }$, and $\dd_\vx\phi$ and $
\dd_\vs \phi $ are defined similarly. We let $\dd\phi =
\paren{\dd_{\vy} \phi, \dd_{\vx} \phi,\dd_\vs \phi}$.

% \[\begin{gathered}\dd_{\vy} \phi = \paren{\partyf{1}{\phi},
% \partyf{2}{\phi}, \cdots, \partyf{n}{\phi} },\ \dd_\vs \phi =
% \paren{\partsif{1}{\phi},\partsif{2}{\phi}, \cdots, \partsif{N}{\phi}
% }\\ \text{ and }\ \dd\phi(\vx,\vs) = \paren{\dd_{\vy} \phi(\vy,
% \vx,\vs), \dd_{\vx} \phi(\vy,\vx,\vs),\dd_\vs
% \phi(\vy,\vx,\vs)}\in \rn\times \rn\times\rr^N.\end{gathered}\]

We use the convenient notation that if $A\subset \rr^m$, then $\dot{A}
= A\smo$.

The singularities of a function and the directions in which they occur
are described by the wavefront set \cite[page
16]{duistermaat1996fourier}: 
\begin{definition}
\label{WF} Let $X$ Let an open subset of $\rn$ and let $f$ be a
distribution in $\mathcal{D}'(X)$.  Let $(\vx_0,\vxi_0)\in X\times
\drn$.  Then $f$ is \emph{smooth at $\vx_0$ in direction $\vxio$} if
there exists a neighbourhood $U$ of $\vx_0$ and $V$ of $\vxi_0$ such
that for every $\phi\in \Dc(U)$ and $N\in\mathbb{R}$ there exists a
constant $C_N$ such that for all $\vxi\in V$,
\begin{equation}
\left|\Fc(\phi f)(\lambda\vxi)\right|\leq C_N(1+\abs{\lambda})^{-N}.
\end{equation}
The pair $(\vx_0,\vxio)$ is in the \emph{wavefront set,} $\wf(f)$, if
$f$ is not smooth at $\vx_0$ in direction $\vxio$.
\end{definition}
 This definition follows the intuitive idea that the elements of
$\WF(f)$ are the point--normal vector pairs above points of $X$ at
which $f$ has singularities.  For example, if $f$ is the
characteristic function of the unit ball in $\mathbb{R}^3$, then its
wavefront set is $\WF(f)=\{(\vx,t\vx): \vx\in S^{2}, t\neq 0\}$, the
set of points on a sphere paired with the corresponding normal vectors
to the sphere.

%\begin{equation}

%\end{equation}
%That is, 

The wavefront set of a distribution on $X$ is normally defined as a
subset the cotangent bundle $T^*(X)$ so it is invariant under
diffeomorphisms, but we do not need this invariance, so we will
continue to identify $T^*(X) = X \times \rn$ and consider $\WF(f)$ as
a subset of $X\times \drn$.

%Let $X$ and $Y$ be open subsets of $\rn$, $m \in\mathbb{R}$.

 \begin{definition}[{\cite[Definition 7.8.1]{hormanderI}}] We define
 $S^m(Y\times X\times \mathbb{R}^N)$ to be the
set of $a\in \Ec(Y\times X\times \mathbb{R}^N)$ such that for every
compact set $K\subset Y\times X$ and all multi--indices $\alpha,
\beta, \gamma$ the bound
\[
\left|\partial^{\gamma}_{\vy}\partial^{\beta}_{\vx}\partial^{\alpha}_{\vsig}a(\vy,\vx,\vsig)\right|\leq
C_{K,\alpha,\beta,\gamma}(1+\norm{\vsig})^{m-|\alpha|},\ \ \ (\vy,\vx)\in K,\
\vsig\in\mathbb{R}^N,
\]
holds for some constant $C_{K,\alpha,\beta,\gamma}>0$. 

 The elements of $S^m$ are called \emph{symbols} of order $m$.  Note
that these symbols are sometimes denoted $S^m_{1,0}$.  The symbol
$a\in S^m(Y,X,\rr^N)$ is \emph{elliptic} if for each compact set
$K\subset Y\times X$, there is a $C_K>0$ and $M>0$ such that
\bel{def:elliptic} \abs{a(\vy,\vx,\vsig)}\geq C_K(1+\norm{\vsig})^m,\
\ \ (\vy,\vx)\in K,\ \norm{\vsig}\geq M.\ee 
\end{definition}

\begin{definition}[{\cite[Definition
        21.2.15]{hormanderIII}}] \label{phasedef}
A function $\phi=\phi(\vy,\vx,\vsig)\in
\Ec(Y\times X\times\dot{\mathbb{R}^N})$ is a \emph{phase
function} if $\phi(\vy,\vx,\lambda\vsig)=\lambda\phi(\vy,\vx,\vsig)$, $\forall
\lambda>0$ and $\mathrm{d}\phi$ is nowhere zero. A phase function is
\emph{clean} if the critical set $\Sigma_\phi = \{ (\vy,\vx,\vsig) \ : \
\mathrm{d}_\vsig \phi(\vy,\vx,\vsig) = 0 \}$ is a smooth manifold with tangent
space defined by $\mathrm{d} \paren{\mathrm{d}_\vsig \phi}= 0$.
\end{definition}
\noindent By the implicit function theorem the requirement for a phase
function to be clean is satisfied if
$\mathrm{d}\paren{\mathrm{d}_\vsig
\phi}$ has constant rank.

\begin{definition}[{\cite[Definition 21.2.15]{hormanderIII} and
      \cite[section 25.2]{hormander}}]\label{def:canon} Let $X$ and
$Y$ be open subsets of $\rn$. Let $\phi\in \Ec\paren{Y \times X \times
{\rr}^N}$ be a clean phase function.  In addition, we assume that
$\phi$ is \emph{nondegenerate} in the following sense:
\[\text{$\dd_{\vy,\vsig}\phi$ and $\dd_{\vx,\vsig}\phi$ are never zero.}\]
  The
\emph{critical set of $\phi$} is
\[\Sigma_\phi=\{(\vy,\vx,\vsig)\in Y\times X\times\dot{\mathbb{R}^N}
: \dd_{\vsig}\phi=0\}.\]  The
\emph{canonical relation parametrised by $\phi$} is defined as
\bel{def:Cgenl} \begin{aligned} \Cc=&\sparen{
\paren{\paren{\vy,\dd_{\vy}\phi(\vy,\vx,\vsig)};\paren{\vx,-\dd_{\vx}\phi(\vy,\vx,\vsig)}}:(\vy,\vx,\vsig)\in
\Sigma_{\phi}},
% &\hspace{1.5cm} \vs\in \rr^N\smo,   
\end{aligned}
\end{equation}
\end{definition}

\begin{definition}\label{FIOdef}
Let $X$ and $Y$ be open subsets of $\rn$. A \emph{Fourier integral operator
(FIO)} of order $m + N/2 - n/2$ is an operator $A:\Dc(X)\to
\mathcal{D}'(Y)$ with Schwartz kernel given by an oscillatory integral
of the form
\begin{equation} \label{oscint}
K_A(\vy,\vx)=\int_{\mathbb{R}^N}
e^{i\phi(\vy,\vx,\vsig)}a(\vy,\vx,\vsig) \mathrm{d}\vsig,
\end{equation}
where $\phi$ is a clean nondegenerate phase function and $a$ is a
symbol in $S^m(Y \times X \times \mathbb{R}^N)$. The \emph{canonical
relation of $A$} is the canonical relation of $\phi$ defined in
\eqref{def:Cgenl}.

The FIO $A$ is \emph{elliptic} if its symbol is elliptic.
\end{definition}

This is a simplified version of the definition of FIO in \cite[section
2.4]{duist} or \cite[section 25.2]{hormander} that is suitable for our
purposes since our phase functions are global.  Because we assume
phase functions are nondegenerate, our FIO can be extended from as
maps from $\Dc(X)$ to $\Ec(Y)$ to maps from $\Ec'(X)$ to $\Dc'(Y)$,
and sometimes larger sets.  For general information about FIOs see
\cite{duist, hormander, hormanderIII}.

The composition of sets is defined as follows.  Let $X$ and $Y$ be
sets and let $A\subset X$ and $B\subset Y\times X$ the composition
\[\begin{aligned}B\circ A &= \sparen{y\in Y\st \exists x\in X,\
(y,x)\in B}\\
B^t &= \sparen{(x,y)\st (y,x)\in B}.\end{aligned}\]

The H\"ormander-Sato Lemma  provides the relationship between the
wavefront set of distributions and their images under FIO.

\begin{theorem}[{\cite[Theorem 8.2.13]{hormanderI}}]\label{thm:HS} Let $f\in \Ec'(X)$ and
let $F:\Ec'(X)\to \Dc'(Y)$ be an FIO with canonical relation $\Cc$.
Then, $\wf(Ff)\subset \Cc\circ \wf(f)$.\end{theorem}

\begin{definition}
\label{defproj} Let $\Cc\subset T^*(Y\times X)$ be the canonical
relation associated to the FIO $A:\mathcal{E}'(X)\to \mathcal{D}'(Y)$.
Let $\Pi_L$ and $\Pi_R$ denote the natural left- and right-projections
of $\Cc$, projecting onto the appropriate coordinates: $\Pi_L:\Cc\to
T^*(Y)$ and $\Pi_R : \Cc\to T^*(X)$.
\end{definition}

Because $\phi$ is nondegenerate, the projections do not map to the
zero section.  
% 
% We have the following result from \cite{hormander}.
% \begin{proposition}
% \label{prop1}
% Let $\dim(X)=\dim(Y)$. Then at any point in $\Cc$:
% \begin{enumerate}[(i)]
% \item if one of $\Pi_L$ or $\Pi_R$ is a local diffeomorphism, then the
% other map is a local diffeomorphism (so $\Cc$ is a local canonical
% graph); 
% 
% \item if one of the projections $\Pi_R$ or $\Pi_L$ is singular at a
% point in $\Cc$, then so is the other. The type of the singularity may
% be different but both projections drop rank on the same set
% \begin{equation}
% \Sigma=\{(\vy,\eta; \vx,\vsig)\in \Cc :
% \det(\mathrm{d}\Pi_L)=0\}=\{(\vy,\eta; \vx,\vsig)\in \Cc : \det
% (\mathrm{d}\Pi_R)=0\}.
% \end{equation}
% \end{enumerate}
% \end{proposition}

If a FIO $\Fc$ satisfies our next definition, then $\Fc^* \Fc$ (or, if
$\Fc$ does not map to $\Ec'(Y)$, then $\Fc^* \psi \Fc$ for an
appropriate cutoff $\psi$) is a pseudodifferential operator
\cite{GS1977, quinto}.

\begin{definition}\label{def:bolker} Let
$\Fc:\Ec'(X)\to \Dc'(Y)$ be a FIO with canonical relation $\Cc$ then
$\Fc$ (or $\Cc$) satisfies the \emph{semi-global Bolker Assumption} if
the natural projection $\Pi_Y:\Cc\to T^*(Y)$ is an embedding
(injective immersion).\end{definition}

%% file: MainThm.tex
\section{The Main Theorem}\label{sect:main}

In this section we define our transform and give conditions under
which our transform satisfies the Bolker Assumption.  We consider a Radon transform in $\rn$ that is a generalization of the
transforms studied in \cite{vline,Terzioglu-KK-cone,Web,
zhang2020recovery}. This transform will integrate on surfaces of rotation with vertex on
the $\xn=0$ hyperplane. 

We start with a function that will  define the
surfaces.  \bel{hyp1}\begin{gathered} \text{Let
$q:[0,\infty)\to[0,\infty)$ be continuous on $[0,\infty)$, $C^\infty$
on $\oinf$, }\\
\text{$q(0)=0$, and $q(r)>0$ if $r>0$.}\end{gathered}\ee

% \jc{\tred{I get now why we can use $q(r)>0$ for $r>0$ since you
% \define $g>0$ later. $q>0$ is more general than $q'>0$ so I like
% \yours more.}}

If $\vx=(x_1,x_2,\dots,x_n)\in \rn$ then we let $\vx' =
(x_1,x_2,\dots, x_{n-1})\in \rnmo$ so $\vx = (\vx',\xn)$.  Now, let
$X=\sparen{(\vx',\xn)\st \vx'\in \rnmo,\xn\in (0,\infty)}$ denote the
half-space $\xn>0$ in $\rn$.  Let $Y=(0,\infty)\times \rnmo$.  Then,
for $(E,\vxo)\in Y$, the surface of integration of our Radon transform
is given by \bel{def:S}S(E,\vxo)= \sparen{(\vx',\xn)\st
\xn=Eq\paren{\norm{\vx'-\vxo}}\ \vx'\in \rnmo\setminus
\sparen{\vxo}}.\ee Note that $S(E,\vxo)$ has axis of rotation
$\{(\vxo,\xn)\st \xn>0\}$ and vertex at $(\vxo,0)$ (which is not in
$S(E,\vxo)$).  The surface $S(E,\vxo)$ is characterized by the
equation \bel{def:Psi}\begin{gathered} \Psi(E,\vxo,(\vx',\xn))=0
\text{\ \ where}\\
\Psi(E,\vxo,(\vx',\xn)): = \xn-Eq\paren{\norm{\vx'-\vxo}}
\end{gathered}\ee

The \emph{generalized cone Radon transform} is given, for $f\in
L^2_c(X)$, by \bel{def:R}
\begin{aligned}Rf(E,\vxo)&=\int_{\vx\in S(E,\vxo)} f(\vx)\dd S(\vx)\\
& =\int_{(\vx,\xn)\in X} f(\vx',\xn)\,\norm{\nabla_\vx
\Psi}\,\delta\paren{\Psi(E,\vxo,(\vx',\xn))}\dd \vx'\,\dd
\xn\end{aligned}\ee where we use \cite[eq.\
(1)]{palamodov2012uniform} and the relation of the transform $M_\Psi$
in that article to $R$ (see also \cite[\S 6.1]{hormanderI}). Thus,
$Rf(E,\vxo)$ integrates $f$ over the surface of rotation $S(E,\vxo)$
in surface area measure. 

Our first main theorem allows us to analyze mapping properties of $R$
microlocally and in Sobolev space.

% \begin{theorem}\label{thm:BolkerRn}
%   Let $q:[0,\infty) \to [0,\infty)$ satisfy \eqref{hyp1}.  Then, the
% associated Radon transform $R$ is a FIO of order $(n-1)/2$ satisfying
% the Bolker assumption if and only if $g=q'/q$ satisfies 
% \begin{enumerate}[(a)]
% \item \label{g} $g:\oinf\to\oinf$ is injective, and
% \item\label{g'} $g:\oinf\to\oinf$ is an immersion (i.e., $g'(r)\neq
% 0\ \forall r\in \oinf$).
% \end{enumerate}
% \end{theorem}
% 
%  \begin{remark} Let $\Cc$ be the canonical relation of $R$.  Our proof
% will show that condition \ref{g} is equivalent to $\Pi_L:\Cc\to
% T^*(Y)$ being injective.  
% Condition \ref{g'} is equivalent to
% $\Pi_L:\Cc\to T^*(Y)$ being an immersion.  
% 
% Note that condition \ref{g'} implies condition \ref{g}, but by
% including \ref{g}, the theorem becomes an equivalence.
% 
% \end{remark}

% \tc{I'm torn about the theorem.  If people want to check Bolker then
% they need only \ref{g'} and really Bolker is equivalent to \ref{g'}.
% So I rewrote the theorem so you need to check only \ref{g'}.  Then, I
% expand this in the remark.  The original formulations are commented
% above.}

\begin{theorem}\label{thm:BolkerRn}
  Let $q:[0,\infty) \to [0,\infty)$ satisfy \eqref{hyp1}.  Then, the associated generalized cone Radon transform $R$
is an elliptic FIO of order $\frac{1-n}{2}$. 

 Let $g=q'/q$ and assume $g:\oinf\to\oinf$. Then, the transform $R$
satisfies the Bolker assumption if and only if \bel{g'}
\text{$g'(r)\neq 0$ for all $r\in \oinf$ (i.e., $g$ is an immersion).}
\end{equation}
This condition is equivalent to $qq''-(q')^2$ being nowhere zero for
$r\in \oinf$.\end{theorem}

\begin{remark}\label{rem:qq'} First, note that our theorems 
are valid for all Radon transforms defined on surfaces $S(E,\vxo)$ for
which $q$ and $g$ satisfy the conditions in Theorem \ref{thm:BolkerRn}
and the weights on the surfaces are smooth and nowhere zero.  This is
true because our proofs use microlocal analysis, which does not depend
on the specific weight.  Ellipticity of the operator follows because
the weight is assumed to be nowhere zero.

Let $\Cc$ be the canonical relation of $R$.  For the Bolker Assumption
to hold, $\Pi_L:\Cc\to T^*(Y)$ needs to be both injective and
immersive.  In our proof, we will show that \eqref{g'} is equivalent
to $\Pi_L$ being immersive.  We will also show that the condition
\bel{g} g:\oinf\to\oinf\ \, \text{is injective}\ee is equivalent to
$\Pi_L$ being injective.  However, condition \eqref{g'} implies this
new condition \eqref{g} for the following reason: if $g'$ is never $0$
then $g$ must be strictly monotonic because the domain of $g$,
$\oinf$, is connected.
\end{remark}

We will first prove this theorem in $\rtwo$ since this provides the
main ideas. Then, we provide the general proof for $\rn$.

%% file: R2_2.tex
\subsection{Proof of Theorem \ref{thm:BolkerRn} in
$\rtwo$}\label{sect:proofR2}
In this case the ``surface of rotation'' $S(E,\xo)$ consists
of two curves that are mirror images of each other, so we will
introduce two Radon transforms.  Throughout this section we use the
coordinates $(E,x_0,x_1,x_2,\sigma)\in Y\times X \times \dR$.

Define 
\[D_j = D_j(\xo) = \sparen{(x_1,x_2)\st \nojxxo>0, x_2>0}\]
and let $\Psi_j(E,x_0,x_1,x_2)=x_2-Eq(\nojxxo)$. Then, for
$f\in
L^2_c(X)$, we define the transforms $R_j$, for $j=1,2$, as
\begin{equation}
\label{R2Rad}\begin{aligned}R_j f(E,\xo) &= \int_{D_j(\xo)}
\norm{\nabla_{\vx}\Psi_j(E,x_0,\vx)}\delta(x_2-Eq(\noj(x_1-\xo))f(\vx)\dd \vx \\
&=\int_{-\infty}^{\infty} \int_{D_j(\xo)}      
a(\vx,E,x_0)e^{\Phi_j(E,\xo,\vx,\sigma)}
f(\vx)\dd \vx \dd \sigma,
\end{aligned}
\end{equation} where
\[ \Phi_j(E,\xo,\vx,\sigma) = \sigma
\paren{x_2-Eq(\noj(x_1-\xo))},\] 
and
\begin{equation}
\begin{split}
a(\vx,E,x_0)&=\norm{\nabla_{\vx}\Psi_j(E,x_0,\vx)}\\
&=\sqrt{1+E^2(q'((-1)^j(x_1-x_0)))^2}.
\end{split}
\end{equation}
To get the second line of \eqref{R2Rad} we use the Fourier representation of the delta function. Then, in $\rtwo$, the Radon transform $R$ of
\eqref{def:R} can be written
\begin{equation}
Rf(E,x_0)=R_1f(E,x_0)+R_2f(E,x_0).
\end{equation}
It can be shown that the phase function $\Phi_j$ is non-degenerate
(see definition \ref{phasedef}). The calculation of non-degeneracy is
left to the reader.

The amplitude $a$ is smooth, never zero, and not dependent on the
phase variable $\sigma$. Further the partial
derivatives of $a$, of all orders, are bounded on any compact set.
Therefore $a$ is a symbol order zero. It follows that the $R_j$ and
$R=R_1+R_2$ are elliptic FIO order
$O(R),O(R_j)=0+\frac{1}{2}-\frac{2}{2}=-\frac{1}{2}$, using the
formula of Definition \ref{FIOdef}.

Let $\Dc_j = \oinf\times \sparen{(x_0,x_1)\in\rtwo\st
(-1)^j(x_1-\xo)>0}\times \dR$. Then the canonical relations of the $R_j$ are
\[\begin{aligned}
\Cc_j = &\Big\{\big(\overbrace{(E,\xo)}^{\vy},\overbrace{-\sigma q(\noj(x_1-\xo)), \noj\sigma
Eq'(\noj(x_1-\xo))}^{\mathrm{d}_{\vy}\Phi_j};  \\
&\qquad \vx,\underbrace{\noj \sigma Eq'(\nojxxo), -\sigma}_{-\mathrm{d}_{\vx}\Phi_j}
 \big) \st (E,\xo,x_1,\sigma)\in \mathcal{D}_j, x_2=Eq(\nojxxo)\Big\}\end{aligned}\] 
In these coordinates using $\mathcal{D}_j$, the left projection $\Pi^{(j)}_L : \mathcal{D}_j\to \Pi^{(j)}_L(\mathcal{D}_j)$ of $R_j$ is \[\pilj(E,\xo,x_1,\sigma)=
(E,\xo,-\sigma q(\noj(x_1-\xo)), \noj\sigma Eq'(\noj(x_1-\xo))).\]
Then the left projection $\Pi_L : \mathcal{D}_1\cup\mathcal{D}_2\to \Pi_L\paren{\mathcal{D}_1\cup\mathcal{D}_2}$ of $R$ is defined by $\Pi_L=\Pi^{(1)}_L$ on $\mathcal{D}_1$, and $\Pi_L=\Pi^{(2)}_L$ on $\mathcal{D}_2$. The canonical relation of $R$ is the disjoint union $\Cc=\Cc_1\cup\Cc_2$.

We will now show that condition \ref{g'} is equivalent to $\Pi_L$ an immersion. To do this we consider the derivatives of the $\Pi^{(j)}_L$,
\begin{equation}\label{DPiL}
D\Pi^{(j)}_L=\begin{pmatrix}
  1 & 0 & 0&0\\
  0 & 1 & 0&0\\
  a_{3,1}  & a_{3,2}  &  \nojp\sigma q'(\nojxxo)&-q(\nojxxo)\\
 a_{4,1} & a_{4,2} & \sigma Eq''(\nojxxo)&\noj Eq'(\nojxxo)
\end{pmatrix}.
\end{equation}
The determinant is
\begin{equation}
\begin{split}
\text{det}D\Pi^{(j)}_L&=\text{det}\bpm \nojp\sigma q'(\nojxxo)&-q(\nojxxo)\\
\sigma Eq''(\nojxxo)&\noj Eq'(\nojxxo)
\epm\\
&=\sigma E\paren{q(\nojxxo)q''(\nojxxo)-q'(\nojxxo)^2},
\end{split}
\end{equation}
which is non-vanishing if and only if
\begin{equation}\label{q q' q''}
q(x_1)q''(x_1)-q'(x_1)^2\neq 0,\ \ \ \forall x_1\in\dR.
\end{equation}
Now
$$g'(x_1)=\frac{q''(x_1)}{q(x_1)}-\frac{q'(x_1)^2}{q^2(x_1)}=0\iff q(x_1)q''(x_1)-q'(x_1)^2= 0$$
for $x\in\dR$. The results follows. 

We now show that condition \ref{g} is equivalent to $\Pi_L$ injective. We first consider the implication $\ref{g}\implies \Pi_L\ \ \text{injective}$. 

Let $g$ be injective, and let $(E_1,x_0,x_1,\sigma_1),
(E_2,x'_0,x'_1,\sigma_2)\in\mathcal{D}_j$ be such that\hfil\newline
$\Pi^{(j)}_L(E_1,x_0,x_1,\sigma_1)=\Pi^{(j)}_L(E_2,x'_0,x'_1,\sigma_2)$. Then $E_1=E_2=E$, $x_0=x'_0$, and
\begin{equation}
\label{piL1}
\begin{pmatrix}
-\sigma_1 q((-1)^j(x_1-x_0))\\
(-1)^j\sigma_1 Eq'((-1)^j(x_1-x_0))\\
\end{pmatrix}=\begin{pmatrix}
-\sigma_2 q((-1)^j(x'_1-x_0))\\
(-1)^j\sigma_2 Eq'((-1)^j(x'_1-x_0))\\
\end{pmatrix}.
\end{equation}
It follows that
$$(-1)^{j-1}Eg((-1)^j(x_1-x_0))=(-1)^{j-1}E g((-1)^j(x'_1-x_0)).$$
Hence $x_1=x'_1$, for $j=1,2$, since $E> 0$ and $g$ is injective. Now $\sigma_1 q((-1)^j(x_1-x_0))=\sigma_2 q((-1)^j(x_1-x_0))\implies \sigma_1=\sigma_2$ since $q((-1)^j(x_1-x_0))>0$ on $\mathcal{D}_j$. Thus $\Pi^{(j)}_L$ is injective, for $j=1,2$.

Now let $(E_1,x_0,x_1,\sigma_1)\in \mathcal{D}_1$ and $(E_2,x'_0,x'_1,\sigma_2)\in \mathcal{D}_2$ be such that
$\Pi^{(1)}_L(E_1,x_0,x_1,\sigma_1)=\Pi^{(2)}_L(E_2,x'_0,x'_1,\sigma_2)$. Then $E_1=E_2=E$, $x_0=x'_0$, and
\begin{equation}\label{C1 C2}
\begin{pmatrix}
-\sigma_1 q(-(x_1-x_0)))\\
 -\sigma_1 Eq'(-(x_1-x_0))
\end{pmatrix}=\begin{pmatrix}
-\sigma_2 q(x'_1-x_0)\\
 \sigma_2 Eq'(x'_1-x_0)
\end{pmatrix}.
\end{equation}
Thus it follows that
\bel{negative}g(-(x_1-x_0))=-g(x'_1-x_0).\ee
Now $x_1-x_0<0$ on $\mathcal{D}_1$ and $x'_1-x_0>0$ on
$\mathcal{D}_2$. Further $g(r)>0$ for all $r>0$ by assumption, so
\eqref{negative} is impossible.
Hence $\Pi_L$ is injective on $\Cc$.

%$g(-x_1)=\frac{q'(-x_1)}{q(-x_1)}>0$ for
%$x_1<0$ and $-g(x_1)<0$ for $x_1>0$, since $g>0$ by \ref{g}.  Thus
%$g(-(x_1-x_0))>0$ and $-g(x'_1-x_0)<0$, and we have a contradiction.

We now prove the converse implication, namely $\Pi_L\ \ \text{injective}\implies \ref{g}$. Let $g$ be non-injective, and let $r_1,r_2\in(0,\infty)$ be such that $g(r_1)=g(r_2)$, with $r_1\neq r_2$. We have $\Pi^{(j)}_L(E,x_0,x_1,\sigma_1)=\Pi^{(j)}_L(E,x_0,x'_1,\sigma_2)\iff \ref{piL1}\ \ \text{holds}$. We can write the equations \ref{piL1} as
$$A\boldsymbol{\sigma}=\begin{pmatrix} q((-1)^j(x_1-x_0)) & -q((-1)^j(x'_1-x_0))\\
(-1)^{j-1} q'((-1)^j(x_1-x_0)) & (-1)^j q'((-1)^j(x'_1-x_0))
\end{pmatrix}\begin{pmatrix}\sigma_1\\
\sigma_2
\end{pmatrix}=\begin{pmatrix}0\\
0
\end{pmatrix}.$$
The determinant of $A$ is 
$$\text{det}(A)=(-1)^j \paren{q'((-1)^j(x'_1-x_0))q((-1)^j(x_1-x_0)) -q((-1)^j(x'_1-x_0))q'((-1)^j(x_1-x_0))}.$$
Thus setting $x_1=x_0+(-1)^jr_1$ and $x'_1=x_0+(-1)^jr_2$ (note $x_1\neq x'_1$) yields $\text{det}(A)=0$, since
$$g(r_1)=g(r_2)\implies \frac{q'(r_1)}{q(r_1)}=\frac{q'(r_2)}{q(r_2)}\implies q'(r_1)q(r_2)=q'(r_2)q(r_1).$$
Hence there exist $\sigma_1, \sigma_2\neq 0$, such that $\boldsymbol{\sigma}\in\text{null}(A)$. For example, $\sigma_1=1$ and
$$\sigma_2=\frac{\sqrt{q^2(r_1)+q'(r_1)^2}}{\sqrt{q^2(r_2)+q'(r_2)^2}}\neq 0$$
is sufficient. Therefore $\Pi_L$ is non-injective. Finally, $\ref{g'} \implies \ref{g}$ (see Remark \ref{rem:qq'}), so condition \ref{g'} is equivalent to the Bolker
Assumption. This completes the proof.  $\qed$

\begin{remark}
\label{rem3.3} When $g$ is non-injective, $\Pi_L$ is non-injective
(see Remark \ref{rem:qq'}) and artifacts can be generated.  Using
\eqref{def:C} one can show that $\Cc^t\circ \Cc \subset
\Delta \cup \Lambda $ where $\Delta$ is the diagonal in $T^*X\times
T^*X$ and $\Lambda\subset T^*X\times T^*X$.  This is important because,
if $\lambda\in \wf(f)$ then
\[\Cc^t\circ\Cc\circ\sparen{\lambda}\subset\paren{\Delta\circ
\sparen{\lambda}}\cup
\paren{\Lambda\circ\sparen{\lambda}}=\sparen{\lambda}\cup
\Lambda\circ\sparen{\lambda},\] and by the H\"ormander-Sato Lemma both
a visible singularity at $\lambda$ and an artifact at
$\Lambda\circ{\sparen{\lambda}}$ could be in $\wf(R^*\psi Rf)$ (where
$\psi$ is a cutoff to make $R^*\psi R$ defined).

To describe the artifacts and, implicitly, $\Lambda$, note that
$\Cc_j^t\circ\Cc_i=\emptyset$ for $i\neq j$ (see \eqref{C1 C2} and
\eqref{negative}). This means that
$\Cc^t\circ\Cc=(\Cc_1\cup\Cc_2)^t\circ(\Cc_1\cup\Cc_2)=
\paren{\Cc_1^t\circ\Cc_1}\cup \paren{\Cc_2^t\circ\Cc_2}.$ Therefore,
any artefacts are due to the $\Cc_j^t\circ\Cc_j
\subset\Delta\cup\Lambda$ for $j=1,2$.

Since we now assume $g$ is not injective, we can choose $r_1\neq r_2$
such that $g(r_1) = g(r_2)$.  Let $f$ be a distribution and assume for
some $(E,\xo,\sigma_1)$,
$$\lambda = \paren{(\xo + \noj r_1, Eq((-1)^jr_1)),\sigma_1 Eq'(\noj r_1), -\sigma_1}\in
\text{WF}(f)
$$Equivalently assume there exists an integration curve which intersects
a singularity of $f$ normal to its direction.  Then, one can choose a
$\sigma_2\in \dR$ such that
$\Pi_L(E,\xo,r_1,\sigma_1)=\Pi_L(E,\xo,r_2,\sigma_2)$ where we are
using the coordinates above \eqref{piL1}.  This is true because the
ratio $g=q'/q$ is the same at $r_1$ and $r_2$.  Then by calculating
$\Cc_j^t\circ \Cc_j\circ \sparen{\lambda}$ and using the
H\"ormander-Sato Lemma, one sees that $R_j^*R_jf$ could have an
artefact  at $(\xo+\noj r_2, Eq(r_2), \sigma_2\noj
Eq'(r_2),-\sigma_2)$, and this implicitly describes $\Lambda$.  Note
that artefacts can occur for each $R_j$.  These artifacts will be
shown in simulations in Example \ref{ex4}.
\end{remark}

%% file: Rn.tex
\subsection{Proof of Theorem \ref{thm:BolkerRn} in
$\rn$}\label{sect:proofRn}

\begin{proof}
Let $f\in L^2_c(X)$, then the Radon transform $Rf$ in \eqref{def:R}
satisfies \bel{def:R-FIO}
\begin{split}
Rf(E,\vxo)&=\int_{\mathbb{R}^n}\delta(\Psi(E,\vxo,\vx))\,\norm{\nabla_{\vx}
\Psi(E,\vxo,\vx)}f(\vx)
\dd\vx\\
&=\int_{-\infty}^{\infty}\int_{\mathbb{R}^n}e^{-i\sigma\paren{\xn-Eq\paren{\norm{\vx'-\vxo}}}}
\norm{\nabla_{\vx}\Psi(E,\vxo,\vx)} f(\vx)
\mathrm{d}\vx\mathrm{d}\sigma,
\end{split} 
\end{equation}
and $R$ is an elliptic FIO of order $\frac{1-n}{2}$ satisfying
Definition \ref{FIOdef} with nondegenerate phase function
\bel{def:Phi}\Phi\paren{(E,\vxo),(\vx',\xn),\sigma} =\sigma
\Psi(E,\vxo,\vx)= \sigma \paren{\xn-Eq\paren{\norm{\vx'-\vxo}}}\ee and
symbol \[a(E,\vxo,\vx,\sigma) = \norm{\nabla_\vx
\Psi(E,\vxo,(\vx',\xn))}=\sqrt{1+\norm{\nabla_{\vx'}
Eq\paren{\norm{\vx'-\vxo}}}^2}\] for the same reasons as for the
transforms $R_j$ in Section \ref{sect:proofR2}.  One uses the
definition of canonical relation \eqref{def:Cgenl} to show that the
canonical relation of $R$ is \bel{def:C}\begin{aligned} \Cc =
\Big\{\Big(\overbrace{(E,\vxo)}^{\vy},\overbrace{-\sigma q(r), -\sigma
Eq'(r)\vom}^{\dd_\vy\Phi};
\overbrace{(\vxo+r\vom,Eq(r))}^{\vx},\overbrace{\sigma E
q'(r)\vom,-\sigma}^{-\dd _\vx \Phi}\Big)\st \\
(E,\vxo)\in Y, r>0, \vom\in \snmt, \sigma\neq
0\Big\},\end{aligned}\ee 
where we have written $\vx' = \vxo+r\vom$.
Note that, $r>0$ for points in $\Cc$ because $(\vxo+r\vom,\xn)\in S(E,\vxo)$,
and therefore $\xn=Eq(r)>0$ by \eqref{hyp1}.

We now investigate the mapping properties of $\Pi_L$.  Note that
$(E,\vxo,r,\vom,\sigma)\in Y\times (0,\infty)\times \snmt\times
\dot{\rr}$ give coordinates on $\Cc$.  In these coordinates, $\Pi_L$
becomes \bel{piL}\begin{gathered} (E,\vxo,r,\vom,\sigma)\mapsto
(E,\vxo,-\sigma q(r),-\sigma Eq'(r)\vom)=(E,\vxo,\eta,\vxi)\\
\text{where} \ \ \eta = -\sigma q(r)\in \dR,\quad \vxi = -\sigma
Eq'(r)\vom\in \rnmo.\end{gathered} \ee Let $\vlam= (E,\vxo,\eta,\vxi)$
be in the image of $\Pi_L$. Then \eqref{piL} gives $(E,\vxo)$ and we
need to find formulas for $r$, $\vom$, and $\sigma$ in terms of
$\vlam$.  Recall that we have defined $g(r) = q'(r)/q(r)$.  Since, by
assumption, $g:\oinf\to\oinf$ and $q(r)>0$ on $(0,\infty)$, $q'$ is
always positive on $\oinf$. This explains why $\eta$ and $\vxi$ are
nonzero.

Let \[\vw = \frac{1}{E\eta}\vxi= g(r)\vom,\] then $\vw$ is known from
\eqref{piL} as is $\norm{\vw} = g(r)$.

 First, assume $g$ is injective.  Then, $r$ is determined and $q(r)$
is known and so \bel{sigma omega}\sigma=\frac{-\eta}{q(r)}\qquad \vom
= \frac{-1}{\sigma E q'(r)}\vxi\ee are determined from
$(E,\vxo,\eta,\vxi)$ using \eqref{piL}.  Therefore, $\Pi_L$ is
injective.  Next, if $g$ is not injective then for multiple values of
$r$, \eqref{piL} maps to the same point and $\Pi_L$ is not injective.
Therefore $g$ is injective if and only if $\Pi_L$ is.

To prove $\Pi_L$ is an immersion if and only if $g'$ is never zero,
one does a calculation similar to the calculations in
\eqref{DPiL}-\eqref{q q' q''}. The calculation is simplified by
choosing orthonormal coordinates on $S^{n-1}$ at $\omega$. The result
is that $\det{D\Pi_L}\neq 0$ if and only if
$(q'(r))^{n-1}\bparen{q(r)q''(r)-(q'(r)^2}\neq 0$ and this is
expression is nonzero if and only if $g'$ is never zero on $\oinf$.

Finally, as noted in Remark \ref{rem:qq'}, condition \ref{g'} implies
condition \ref{g}, so condition \ref{g'} is equivalent to the Bolker
Assumption.

%Orthogonal coordiantes mean take an orthonormal basis of $\rn$
%including $\vom_1,\vom_2,\dots,\vom_{n-1},\vom$ and $\vom$, and use
%coordinates $(x_1,x_2,\dots,x_{n-1}\mapsto
%\sum_{j=1}^{n-1}x_j\vom_j+\sqrt{1-(x_1^2+\dots+x_{n-1}^2)}\vom$

\end{proof}

%% file: Sobolev.tex
\subsection{Sobolev Smoothness}\label{sect:Sobolev}

In this section we describe the microlocal continuity properties of
$R$.  Then, we analyze visible and invisible features in the
reconstruction. First we introduce Sobolev spaces and Sobolev wavefront sets \cite{Petersen,
Rudin:FA}.

%We first define Sobolev space and Sobolev wavefront

%To describe the continuity and stability properties of $R$ we
%introduce Sobolev spaces and Sobolev wavefront set \cite{Petersen,
%Rudin:FA}.

\begin{definition}\label{def:Halpha}Let $\alpha\in \rr$.
  Then $H^\alpha(\rn)$ is the set of all distributions for which their
Fourier transform is a locally integrable function and such that the
Sobolev norm \bel{Hsnorm}\norm{f}_\alpha = \paren{\int_{\vxi\in \rn}
\abs{\Fc(f)(\vxi)}^2\paren{1+\norm{\vxi}^2}^\alpha\,\dd\vxi}^{1/2}<\infty.\ee
Let $\Omega\subset \rn$.  Then, $H^\alpha(\Omega)$ will be the set of all
distributions in $H^\alpha(\rn)$ that are supported in $\Omega$, and
$H^\alpha_c(\Omega)$ will be all those of compact support in $\Omega$. We define $H^\alpha_\loc(\Omega)$ as the set of all distributions $f$ supported
in $\Omega$ such that for each $\vp\in \Dc(\Omega)$, the product $\vp
f\in H^\alpha(\rn)$.
\end{definition}

\noindent We give $H^\alpha_c(\Omega)$ the topology using the Sobolev norm
(so $H^\alpha_c(\Omega)$ is not closed), and we give
$H^\alpha_\loc(\Omega)$ the topology defined by the seminorms
$\norm{f}_{\alpha,\vp} = \norm{\vp f}_\alpha$ (so
$H^\alpha_\loc(\Omega)$ is metrizable).

\begin{definition}\label{def:cont}
Let $m\in \rr$ and let $\Omega'$ be an open set in $\rn$.  Then, the
linear map $F:H^\alpha_c(\Omega)\to H^{\alpha-m}_\loc(\Omega')$ is
continuous if for each $\vp\in \Dc(\Omega)$ and $\tilde{\vp}\in
\Dc(\Omega')$ the product map $\tilde{\vp}\, F\,\vp$ is continuous
from $H^\alpha(\Omega)$ to $H^{\alpha-m}(\Omega')$ (in Sobolev
norms).\end{definition}

\begin{corollary}
\label{cr1} Let $q$ satisfy \eqref{hyp1} as well as condition
 \eqref{g'} of Theorem \ref{thm:BolkerRn}. Then $R$ is continuous
from $H^\alpha_c(X)$ to $H^{\alpha+(n-1)/2}_\loc (Y)$. 
\end{corollary}

\noindent This indicates that the forward map is stable in Sobolev scale
$\frac{1-n}{2}$.
\begin{proof}
The operator $R$ is an FIO of order $\frac{1-n}{2}$ with immersive
left projection and for each $\tilde{\vp}$ and $\vp$ in Definition
\ref{def:cont}, the operator $\tilde{\vp} R \vp$ is \emph{compactly
supported}.  Therefore \cite[Theorem 4.3.1]{Ho1971} can be used to
check our definition of continuity from $H^\alpha_c$ to
$H^{\alpha+(n-1)/2}_\loc$.
\end{proof}

\noindent We now define the Sobolev wavefront sets \cite{Petersen}.  This will
provide the language to describe the strength of the visible
singularities in Sobolev scale.

\begin{definition}\label{def:WFa} Let $\alpha\in \rr$ and let
$X\subset \rn$.  Let $f\in \Dc'(X)$ and let $(\vx,\vxi)\in X\times
\dot{\rn}$.  Then, $f$ is \emph{(Sobolev) smooth to order $\alpha$ at
$(\vx,\xi)$} if there is a smooth cutoff function $\vp$ at $\vxo$ and
a conic neighborhood $V$ of $\vxi$ such that \bel{WFa}\int_{\veta\in
V} \abs{\Fc(\vp
f)(\veta)}^2(1+\norm{\veta}^2)^\alpha\dd\veta<\infty.\ee If $f$ is not
smooth to order $\alpha$ at $(\vx,\vxi)$, then $(\vx,\vxi)$ is in the
\emph{Sobolev wavefront set} $\wf^\alpha(f)$.
\end{definition}

If $V$ were replaced by $\rn$ in the integral \eqref{WFa} then
boundedness of the integral would mean that $\vp f$ is in $H^\alpha$.
By restricting the integral to be over $V$ we require $\vp f$ to be in
$H^\alpha$ \emph{only} in some conic neighborhood of $\vxi$. This is a
Sobolev equivalent of Definition \ref{WF} for $C^\infty$ wavefront
set: rapid decrease in $V$ of the localized Fourier transform is
replaced with finite Sobolev seminorm in $V$.

Our next theorem gives the precise relationship between Sobolev
singularities of $f$ and those of $Rf$.

\begin{theorem}\label{thm:WFa-R} Assume $q:[0,\infty)\to[0,\infty)$
satisfies \eqref{hyp1} and \eqref{g'}. Let $R$ be the associated
generalized cone Radon transform.  Let $(\vx,\vxi)\in \rn\times
\dot{\rn}$ and assume that $\vxi'\neq \zero$ and $\xin\neq 0$.  Then,
\bel{wf corresp}(\vx,\vxi)\in \wf^\alpha(f)\quad \iff\quad
\big(E,\vxo,-\sigma q(r), -\sigma E q'(r)\vom\big)\in
\wf^{\alpha+(n-1)/2}(Rf)\ee where \bel{Y coords}\begin{array}{rlll}
&\vom= -\frac{\xin}{\abs{\xin}} \frac{\vxi'}{\norm{\vxi'}}, \qquad
&\,\,r=g\inv\paren{\frac{\norm{\vxi'}}{\xn\abs{\xi_n}}}, \qquad\sigma
= -\xin, \qquad \ \text{then}\\
&E=\frac{\xn}{q(r)},&\vxo=\vx'-r\vom &\phantom{x}
\end{array}\ee
\end{theorem}

In general, Radon transforms smooth singularities (see \cite[Theorem 3.1]{quinto}). Theorem \ref{thm:WFa-R} shows that every singularity of $f$ generates a singularity of
$Rf$ in a specific wavefront direction that is $(n-1)/2$ degrees
smoother in Sobolev scale. Every $(\vx,\vxi)\in \wf^\alpha(f)$ with
$\vxi'\neq \zero$ and $\xin\neq 0$ will create a specific singularity
in $\wf^{\alpha+(n-1)/2}(Rf)$ given by \eqref{Y coords}.  Our proof,
in particular \eqref{pir}, will show that $\Pi_R$ is injective,

\begin{remark}\label{rem:visible wf}  
Vertical and horizontal covectors $(\vx, \vxi)\in \wf(f)$ (where
$\vxi'=\zero$, respectively $\xin=0$) will not create singularities in
$Rf$.  

The reason is as follows.  For a singularity $(\vx,\vxi)$ to be
visible in $Rf$, it must be in the image $\Pi_R(\Cc)$ because
\[\wf(Rf) = \Cc\circ \wf(f)=\Pi_L\circ \Pi_R\inv(\wf(f))\] by
ellipticity, the Bolker assumption, and the H\"ormander-Sato Lemma,
Theorem \ref{thm:HS}.  For $(\vx,\vxi)$ to be in the image of $\Pi_R$,
$\vxi'=\sigma Eq'(r)\vom$ must be nonzero, and this explains why no
vertical covector is in the image of $\Pi_R$.  Furthermore, $\xin =
-\sigma$ must be nonzero, and this explains why no horizontal covector
generates a singularity in $Rf$.  

Therefore, one would expect that those singularities, such as vertical
or horizontal object boundaries, would be difficult  to image in
reconstruction methods. For filtered backprojection reconstruction
methods, this follows from the proof of Theorem \ref{thm:lambda}.
\end{remark}

\begin{proof} 
By the H\"ormander-Sato Lemma (Theorem \ref{thm:HS}), $\wf(Rf)\subset
\Cc\circ \wf(f)$.  If $R$ is elliptic and satisfies the Bolker
assumption equality holds: the $\supset$ containment follows from the
H\"ormander-Sato Lemma applied microlocally to microlocally elliptic
parametrices to $R$.  

The Sobolev version of this wavefront equality follows from Sobolev
continuity of $R$ and of this microlocal parametrix; the proof is
given in \cite{Q1993sing} and \cite[Proposition
A.6]{borg2018analyzing} for the classical Radon transform.  That proof
just uses Sobolev continuity of the classical transform, so it can be
adapted with essentially the same arguments to our case with the
Sobolev continuity order of $\frac{1-n}{2}$ (see \cite[Corollary
6.6]{Petersen} for pseudodifferential operators). This allows us to
say that
\[\Cc\circ \paren{\wf^\alpha(f)\cap
\Pi_R(\Cc)}=\wf^{\alpha+(n-1)/2}(Rf).\] 

To finish the proof of Theorem \ref{thm:WFa-R}, we analyze $\Pi_R$ in
the coordinates $(E,\vxo,r,\vom,\sigma)\in Y\times \oinf\times
\snmt\times \dot{\rr}$ used in Theorem \ref{thm:BolkerRn}.  In these
coordinates the map $\Pi_R$ is described by
\bel{pir}(E,\vxo,r,\vom,\sigma)\mapsto
\paren{(\vxo+r\vom,Eq(r)),\sigma Eq'(r)\vom, -\sigma }=(\vx,\vxi).\ee
If one solves \eqref{pir}, for $(E,\vxo,r,\vom,\sigma)$ one gets
\eqref{Y coords} and therefore \[(E,\vxo,-\sigma q(r),-\sigma
Eq'(r)\vom) = \Pi_L\circ \Pi_R\inv(\vx,\vxi)=\Cc\circ (\vx,\vxi).\]
\end{proof}

% \jc{\tred{swapped $[-c, c]$ and $[a,b]$.}}

% \tc{\tblue{I think it needs to be the other way in this section because
% we are saying $(E,\xo)\in \Ac$ so $[a,b]$ needs to be first.  In the
% next section, you have it $(\xo,E)$ on the graphs. I think that's a
% good way of visualizing the sinograms.  Should we switch in the rest
% of the article or make a comment in the next section?  I think
% switching wouldn't be too difficult.}} 

We will now consider the problem in $\rtwo$ and let $\vx = (x_1,x_2)$
denote a point in $\rtwo$. One of the reconstruction methods in the
next section is a truncated Lambda-filtered back-projection (FBP)
using data for $(E,\xo)$ in a rectangle \bel{region}\Ac:= [a, b]\times
[-c, c] \ee where $0<a<b$ and $0<c$.  Let $\chi_\Ac$ be the
characteristic function of $\Ac$.  The generalized Lambda
reconstruction method we use for functions in $\rtwo$ is
\bel{def:L}\Lc f := R^* \paren{\chi_\Ac\frac{\dd^2} {\dd E^2} Rf}.\ee 

To connect this to reconstructions we need to understand what
singularities of $f$ are visible in its reconstruction from $\Lc$.  An
analysis of added artefacts will be done elsewhere. 

\begin{theorem}[Visible Singularities for $\Lc$ in
$\rtwo$]\label{thm:lambda} Assume $q:[0,\infty)\to[0,\infty)$
satisfies \eqref{hyp1} and \eqref{g'}. Let $R$ be the associated
generalized cone Radon transform in $\rtwo$.  Let $\alpha\in \rr$
and let $\Lc$ be given by \eqref{def:L}.  Let $f\in \Ec'(X)$ and
$(\vx,\vxi)\in \wf^\alpha(f)$.  Then, $(\vx,\vxi)\in
\wf^{\alpha-1}(\Lc f)$ if
\begin{enumerate} \item \label{neq 0}$\xi_1\neq 0$ and $\xi_2\neq 0$
and 
\item\label{open} $\frac{x_2}{q(r)}\in (a,b)$ and $x_1-r\om\in
(-c,c)$ where $r=g\inv\paren{\frac{\abs{\xi_1}}{x_2\abs{\xi_2}}}$ and $\om =
-\frac{\xi_1\xi_2}{\abs{\xi_1\xi_2}}$.
\end{enumerate}

If $(\vx,\vxi)$ does not satisfy condition (\ref{neq 0}), or it does
satisfy (\ref{neq 0}) but $\frac{x_2}{q(r)}\notin [a,b]$, or
$x_1-r\om\notin [-c,c]$, then $\Lc f$ will be smooth at $\vx$ in
direction $\vxi$.
\end{theorem}

\begin{remark}\label{rem:artefacts} 
The reconstructions in section \ref{sect:ex} show the visible
singularities and the invisible singularities predicted by Theorem
\ref{thm:lambda}.  For a visualization of the visible singularities
with broken-ray data (i.e. when $q(r)=r$), see figure \ref{lines}. The
range of $x_0$ and $E$ used is chosen to be consistent with the
simulations conducted in section \ref{sect:ex}. We notice a
greater directional coverage for $\vx$ close to $\{x_2=0\}$, and
conversely for $\vx$ moving away from $\{x_2=0\}$.

\begin{figure}[!h]
\centering
\begin{subfigure}{0.49\textwidth}
\includegraphics[width=0.9\linewidth, height=5.2cm, keepaspectratio]{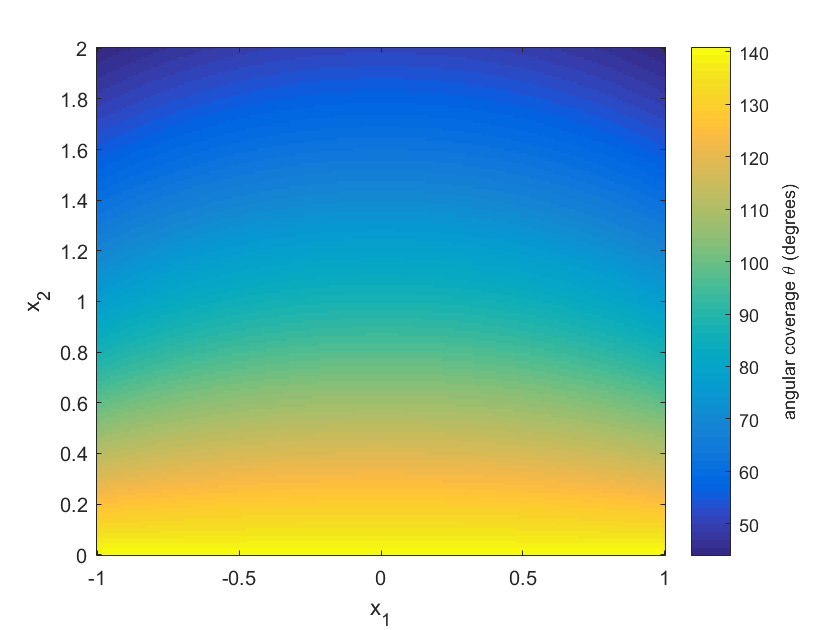}
\subcaption{$\theta$ coverage.}\label{linesA}
\end{subfigure}
\begin{subfigure}{0.49\textwidth}
\includegraphics[width=0.9\linewidth, height=5.2cm, keepaspectratio]{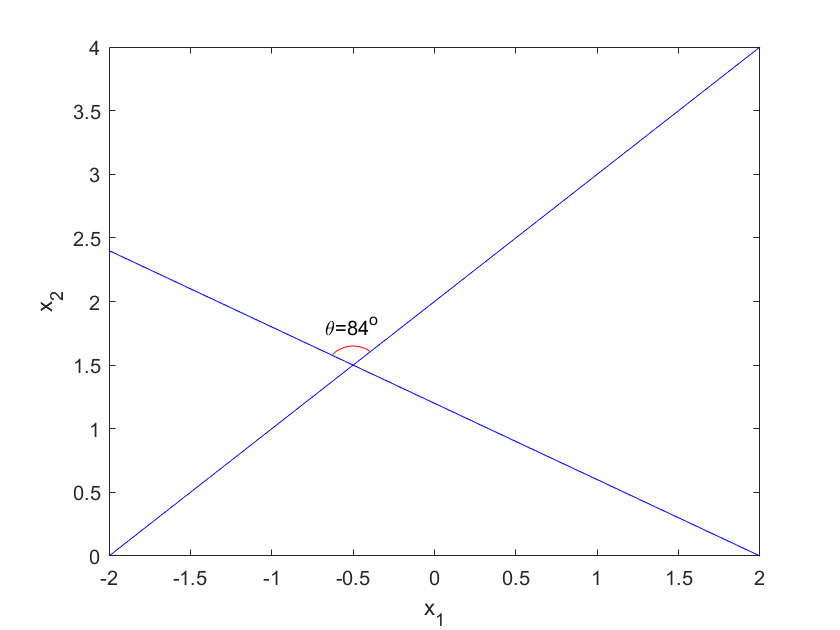}
\subcaption{Example $\theta$.} \label{linesB}
\end{subfigure}
\caption{Left -- set of angular coverage ($\theta$) with broken-ray data on $[-2,2]\times[0,2]$, for $x_0\in[-2,2]$ and $E\in(0,2.83]$ (i.e. when $c=2$, $a=0$ and $b=2.83$). Right -- example $\theta$ for the point $(-0.5,1.5)$. The set of directions resolved by the data ($\xi$) is the red cone displayed in the right hand figure, rotated by $90^{\circ}$ about $(-0.5,1.5)$, minus the direction $\xi=(1,0)$.}
\label{lines}
\end{figure}

Note that this theorem does not
say anything about singularities of $\Lc f$ for $(\vx,\vxi)\in \wf(f)$
for which $\frac{x_2}{q(r)}\in \sparen{a,b}$, or $x_1-r\om = \pm c$.
Singularities at these wavefront directions are more complicated to
analyze and artefacts can be created because of these points.  These
so-called boundary artefacts are seen in figure \ref{FC2}; the
boundary points in the data set labeled ``Edge 1" and ``Edge 2" in figure
\ref{FC2:A} are points in the support of $ R\delta$ at the boundary of
the data set, $\Ac$.  Then, you can see the artefacts they create
(labeled the same way) in figure \ref{FC2:C}.  Similar artefacts are
highlighted in the sinogram in figure \ref{FC2:D} along with the
resulting artefacts in figure \ref{FC2:F}.  These artefacts are
predicted by the microlocal analysis of the operator $\Lc$, and a more
thorough analysis of such artefacts will appear elsewhere. 
\end{remark}

\begin{proof} 
Let $(\vx,\vxi)\in \wf^\alpha(f)$.

First assume that condition \eqref{neq 0} in this theorem holds.
Then, by Theorem \ref{thm:WFa-R}, \[\vlam=(E,\vxo,-\sigma q(r),-\sigma
Eq'(r)\vom)\in \wf^{\alpha+1/2}(Rf)\] by \eqref{Y coords} where
$r=g\inv\paren{\frac{\abs{\xi_1}}{x_2\abs{\xi_2}}}$ and $\om =
-\frac{\xi_1\xi_2}{\abs{\xi_1\xi_2}}$, $E=\frac{x_2}{q(r)}$ and $\xo =
x_1-r\omega$.  Since $\frac{\dd^2}{\dd E^2}$ is elliptic of order two
in this direction, $\vlam\in \wf^{\alpha -3/2}\paren{\frac{\dd^2}{\dd
E^2} Rf}$.

If condition \eqref{open} of this theorem holds, then $(E,\xo)\in
\intt(\Ac)$ and since $\chi_\Ac$ is one in a neighborhood of
$(E,\xo)$, $\vlam\in \wf^{\alpha -3/2}\paren{\chi_\Ac \frac{\dd^2}{\dd
E^2} Rf}$.  Because $R^*$ is elliptic of order $-1/2$ and satisfies
the semi-global Bolker Assumption, $(\vx,\xi) =
\Cc^t\circ\sparen{\vlam}$ is in $\wf^{\alpha-1}(\Lc f)$.  This
statement is proven for $R^*$ using the similar arguments to the
analogous statement for $R$ at the start of the proof of Theorem
\ref{thm:WFa-R}.

Next, if condition \eqref{neq 0} of this theorem holds but
$\frac{x_2}{q(r)}\notin [a,b]$, or $x_1-r\om\notin [-c,c]$, then
$(E,\xo)$ is in the exterior of $\Ac$ and so $\chi_\Ac\frac{\dd^2}{\dd
E^2} Rf$ is zero, hence smooth in a neighborhood of $(E,\xo)$.  Since
$R^*$ is an FIO with canonical relation $\Cc^t$,
$(\vx,\vxi)=\Cc^t\circ{\vlam}\notin \wf(\Lc f)$.

Finally, if \eqref{neq 0} does not hold then, as discussed in Remark
\ref{rem:visible wf} $(\vx,\vxi)$ does not create any singularity in
$Rf$ and therefore not in $\Lc f$ (i.e., $\Cc\circ
\sparen{(\vx,\vxi)}=\emptyset$ and so
$\Cc^t\circ\Cc\circ\sparen{(\vx,\vxi)}=\emptyset$.  This finishes the
proof.
\end{proof}

%% file: ExamplesR2_2.tex
\section{Examples in $\mathbb{R}^2$ and reconstructions}\label{sect:ex} In this section we analyze
several examples to provide perspective on our results. We present reconstructions by inverse crime (noise level zero) to verify our theory. We note in each case if the conditions of Theorem \ref{thm:BolkerRn} (equivalently Bolker) are satisfied.

\begin{example}[CST: Bolker satisfied]
\label{ex1}

Some simple examples of interest are the monomials
\begin{equation}
\label{qex1}
q(r)=r^{\alpha},
\end{equation}
where $\alpha>0$. In this case $g(r)=\frac{\alpha
r^{\alpha-1}}{r^{\alpha}}=\frac{\alpha}{r}$, which is injective on
$(0,\infty)$, and $g'(r)=-\frac{\alpha}{r^2}$, which is never zero.
Hence the conditions of Theorem \ref{thm:BolkerRn} are satisfied and
the Radon integrals satisfy the Bolker assumption.  A specific
monomial of interest in X-ray CT and CST is the straight line
$q_C(r)=r$, when $\alpha=1$. In this case $R_j$, for $j=1,2$, reduces to the well
known line Radon transform in classical X-ray CT, and $R$ reduces to
the broken-ray transforms of
\cite{AmbartsoumianLatifi-Vline2019,vline} in gamma ray source imaging
in CST. See figure \ref{fig1} for example broken-ray integration
curves when $x_0=3$ and $E=0.5$.  Let $B=\{(x,y)\in\mathbb{R}^2 :
x^2+(y-1)^2<0.2^2\}$. Then see figure \ref{FC2} for reconstructions of
a delta function $f=\delta$, centred at $(0,1)$, and disc phantom
$f=\chi_{B}$ from $Rf$, when $q=q_C$.
\begin{figure}[!h]
%\centering
\begin{subfigure}{0.32\textwidth}
\includegraphics[width=0.9\linewidth, height=4cm]{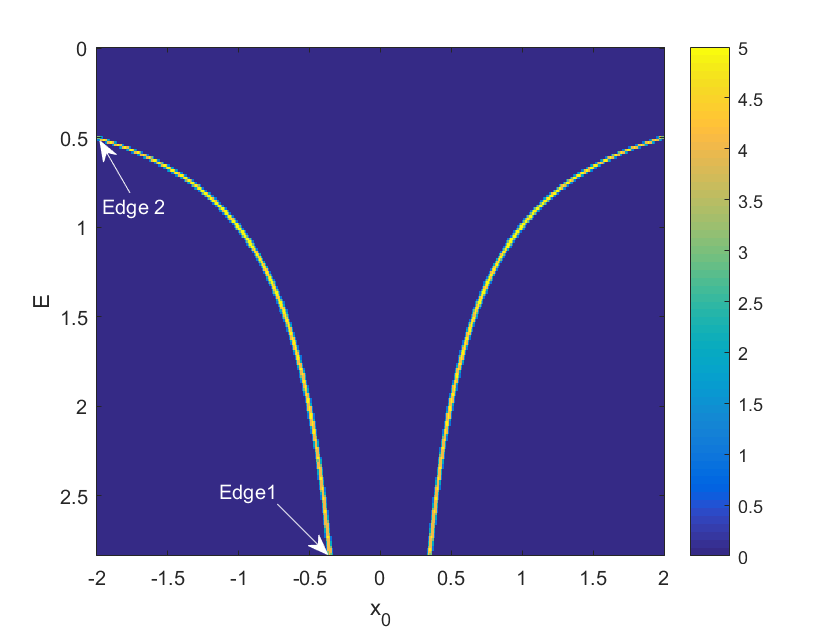}
\subcaption{$R\delta$  sinogram.}\label{FC2:A}
\end{subfigure}
\begin{subfigure}{0.32\textwidth}
\includegraphics[width=0.9\linewidth, height=4cm]{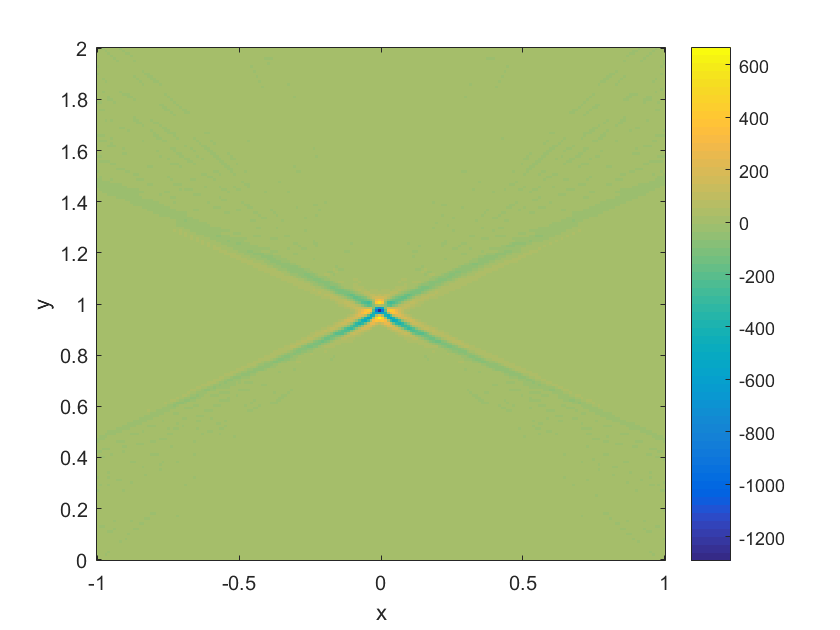}
\subcaption{$R^*\frac{\dd^2}{\dd E^2}R\delta$.}\label{FC2:B}
\end{subfigure}
\begin{subfigure}{0.32\textwidth}
\includegraphics[width=0.9\linewidth, height=4cm]{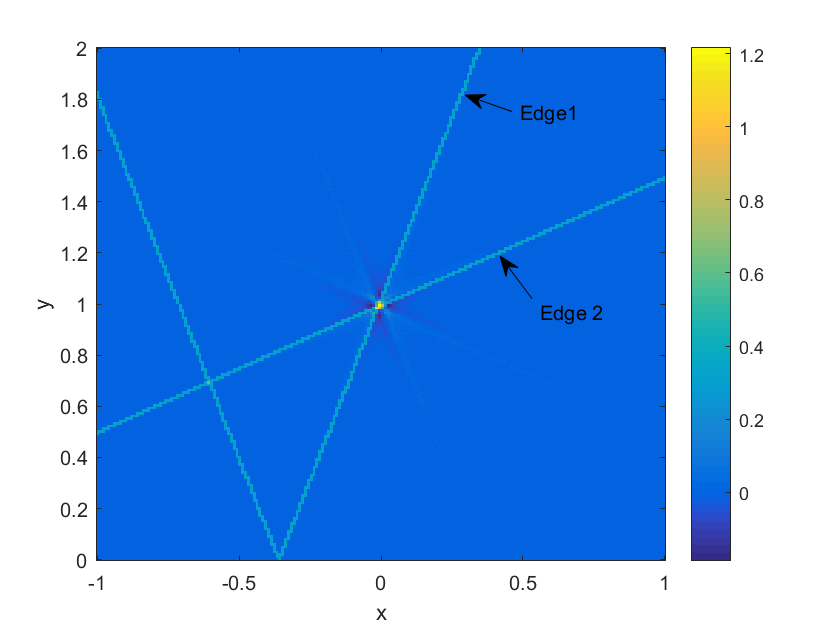}
\subcaption{Landweber.}\label{FC2:C}
\end{subfigure}
\begin{subfigure}{0.32\textwidth}
\includegraphics[width=0.9\linewidth, height=4cm]{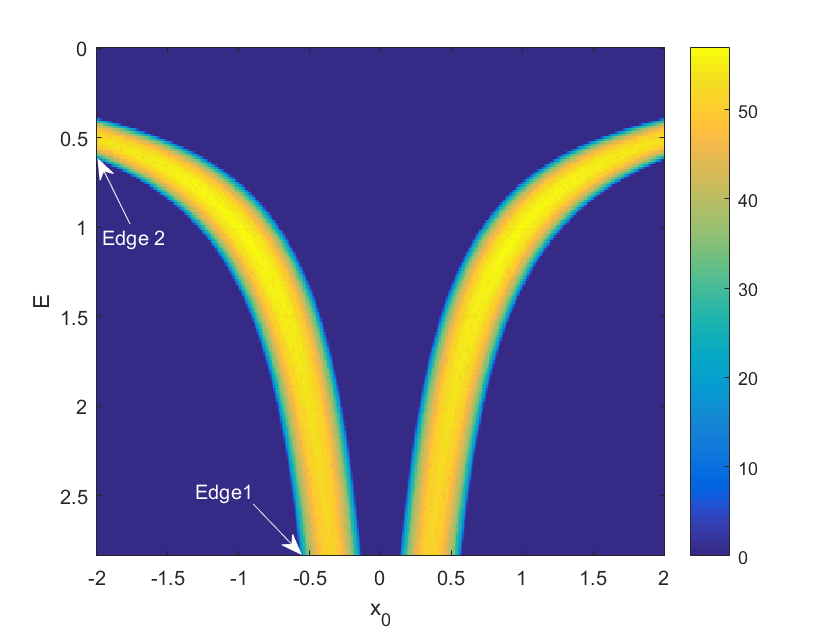}
\subcaption{$R\chi_B$ sinogram.}\label{FC2:D}
\end{subfigure}
\begin{subfigure}{0.32\textwidth}
\includegraphics[width=0.9\linewidth, height=4cm]{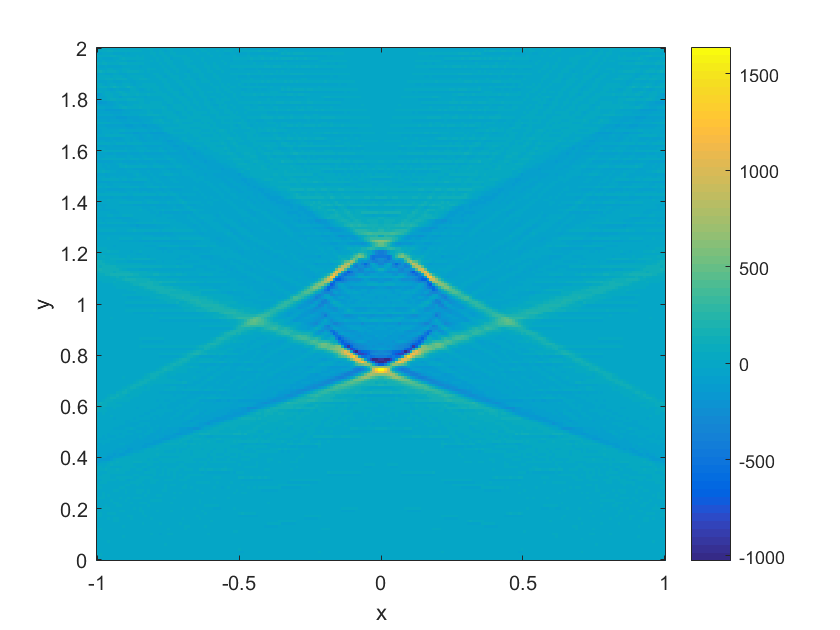}
\subcaption{$R^*\frac{\dd^2}{\dd E^2}R\chi_B$.}\label{FC2:E}
\end{subfigure}
\begin{subfigure}{0.32\textwidth}
\includegraphics[width=0.9\linewidth, height=4cm]{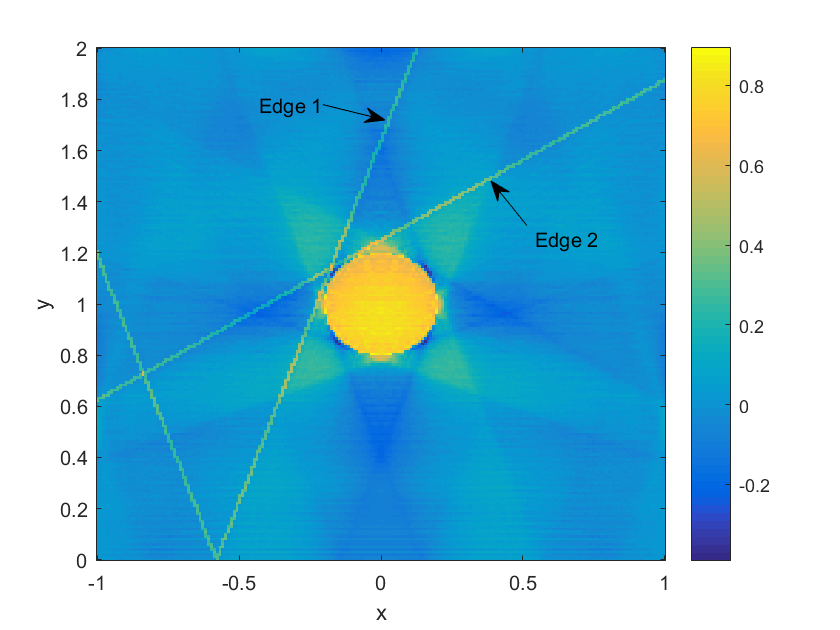}
\subcaption{Landweber.}\label{FC2:F}
\end{subfigure}
\caption{Broken-ray transform reconstructions. Two edges of the data
set are highlighted in each case, which correspond to artefacts in the
reconstruction. These are highlighted in the Landweber image.}
\label{FC2}
\end{figure}
 The scanning region used is $[-1,1]\times [0,2]$, and we simulate
$Rf(E,x_0)$ for $E\in (0,2.83)$ and $x_0\in[-2,2]$. The delta function
is simulated as a characteristic function on a square with small area.
That is $\delta\approx \chi_{S}$, where
$S=[-0.015,0.015]\times[0.985,1.015]$ (i.e. a $3\times 3$ pixel grid
centered on $(0,1)$).  The reconstruction methods used are Filtered
Back-Projection (FBP) (see figures \ref{FC2:B} and \ref{FC2:E} and
remark \ref{rem:artefacts}) using $\frac{\dd^2}{\dd E^2}$ as filter,
and Landweber iteration (see figures \ref{FC2:C} and \ref{FC2:F}).
Note the difference in scales of the color bars between figures
\ref{FC2:B} and \ref{FC2:C}, and \ref{FC2:E} and \ref{FC2:F}. The aim
of the generalized Lambda reconstruction \eqref{def:L} is to
recover the image singularities, and the reconstructed values give
primarily qualitative information. Therefore, these color bar ranges
of figures \ref{FC2:B} and \ref{FC2:E} are chosen to show the
singularities of the object.  The Landweber iteration approximates the
exact solution numerically, and thus the color bar ranges of figures
\ref{FC2:C} and \ref{FC2:F} more closely represent the original density range
%are chosen to correspond to density values
(i.e. $[0,1]$ for the phantoms considered).

We see vertical and horizontal blurring due to limited data in the
Landweber iteration. This is because not all wavefront
directions are visible with broken-ray data, as described by Theorem
\ref{thm:lambda}. This effect is illustrated in figure \ref{linesA}.
We can see that there is only limited angular coverage on the boundary
of $B$ and at $(0,1)$ (the location of $\delta$), where the test
phantoms have singularities. Additionally we see artefacts appearing
along broken-rays at the boundary of the dataset. This is due to the
sharp cutoff in the sinogram space (see figures \ref{FC2:A} and
\ref{FC2:D}). We see similar boundary artefacts occurring in
\cite{borg2018analyzing} in reconstructions from limited line integral
data. Two boundary points in the support of $\chi_A Rf$ are
labelled by ``Edge 1" and ``Edge 2" in the sinograms of figures
\ref{FC2:A} and \ref{FC2:D}, and they generate artifacts that
are shown along broken-ray curves in the image reconstructions of
figures \ref{FC2:C} and \ref{FC2:F}. The broken-ray curves in the
image space are labelled similarly by ``Edge 1" and ``Edge 2", as in
sinogram space.  Note that only half of the broken-ray curve at Edge 2
intersects $[-1,1]\times [0,2]$, and hence the Edge 2 artefacts appear
along lines in figures \ref{FC2:C} and \ref{FC2:F}. Similar boundary
artefacts are observed in \cite{vline}, where the authors present
reconstructions of $\chi_B$ and a Shepp-Logan phantom from $Rf$. The
upper $E$ limit used by \cite{vline} ($E$ is equivalent to the cone
opening angle $\omega$, in the notation of \cite{vline}) is greater
than the maximum $E$ used here however, and hence the reconstructions
\cite{vline} better resolve the horizontal singularities.

If artefacts due to a $\Lambda$ (as in Remark \ref{rem3.3}) are
present, we would expect to see them highlighted in the FBP
reconstruction (as in \cite{webber2020joint}). This is not the case
however and there is no evidence of $\Lambda$ artefacts. This is as
predicted by our theory and is in line with the results of
\cite[Theorem 14]{terzioglu2019some} for a related but
overdetermined transform which show that the normal operator $R^*R$
is a pseudo-differential operator when $q=q_C$
\end{example}

\begin{example}[BST: Bolker satisfied]
\label{ex2}
The curves of integration in BST are \cite{Web}
$$q_B(r)=\frac{r}{\sqrt{r^2+1}}.$$
See figure \ref{fig1} for an example Bragg curve when $x_0=3$ and $E=2$. $q_B$ describes the integration curves for the central scanning profile $x_2=0$, using the notation of \cite{Web}. Explicitly we set $x_2=0$ in \cite[equation (4.2)]{Web} to obtain $q_B$. For an analysis of the general case when $x_2\in(-1,1)$, see appendix \ref{appA}.

The first and second order derivatives of $q_B$ are
$$q'_B(r)=\frac{1}{(r^2+1)^{\frac{3}{2}}}\neq 0,\ \ \ q''_B(r)=-\frac{3r}{(r^2+1)^{\frac{5}{2}}}.$$
Hence $g_B(r)=\frac{q'(r)}{q(r)}=\frac{1}{r(r^2+1)}$, which is injective on $(0,\infty)$, and it follows that
\begin{figure}[!h]
%\centering
\begin{subfigure}{0.32\textwidth}
\includegraphics[width=0.9\linewidth, height=4cm]{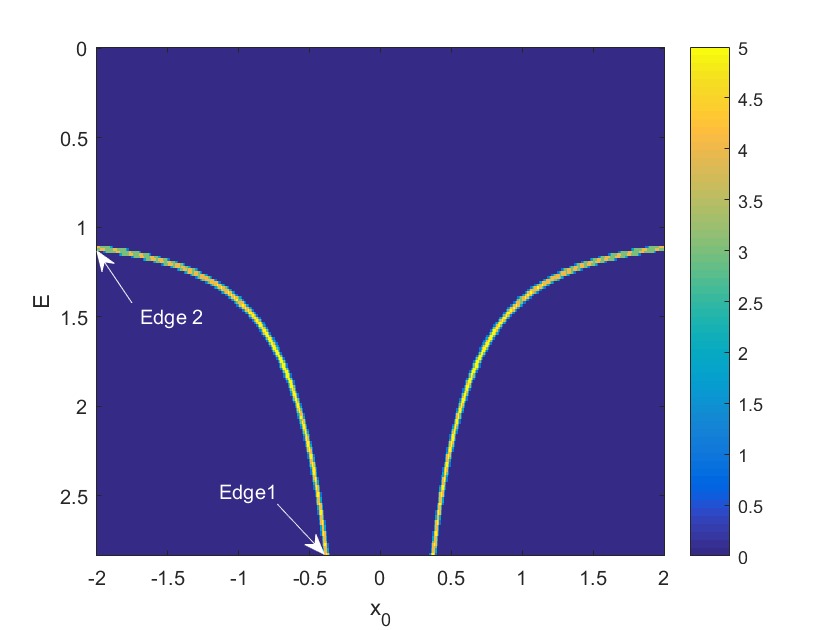}
\subcaption{$R\delta$ sinogram.}\label{FC1:A}
\end{subfigure}
\begin{subfigure}{0.32\textwidth}
\includegraphics[width=0.9\linewidth, height=4cm]{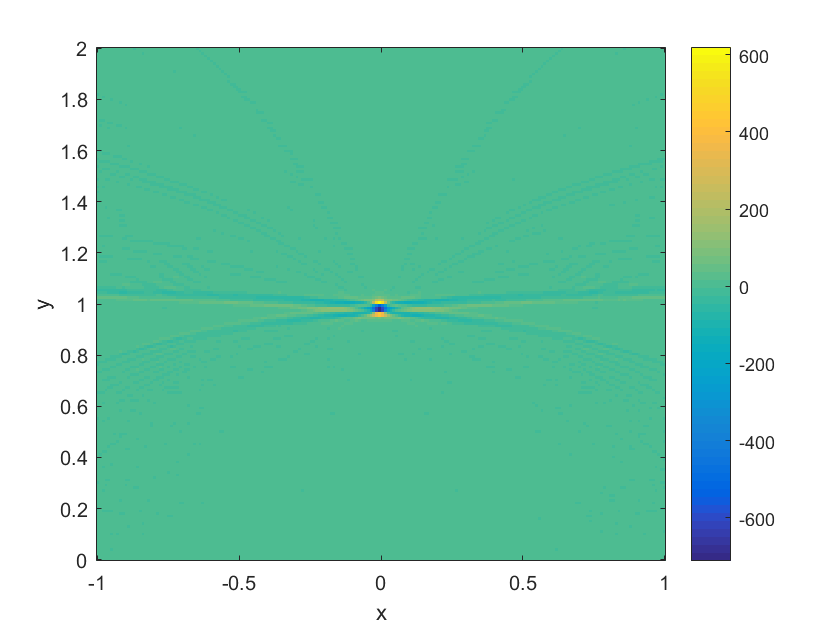}
\subcaption{$R^*\frac{\dd^2}{\dd E^2}R\delta$.}\label{FC1:B}
\end{subfigure}
\begin{subfigure}{0.32\textwidth}
\includegraphics[width=0.9\linewidth, height=4cm]{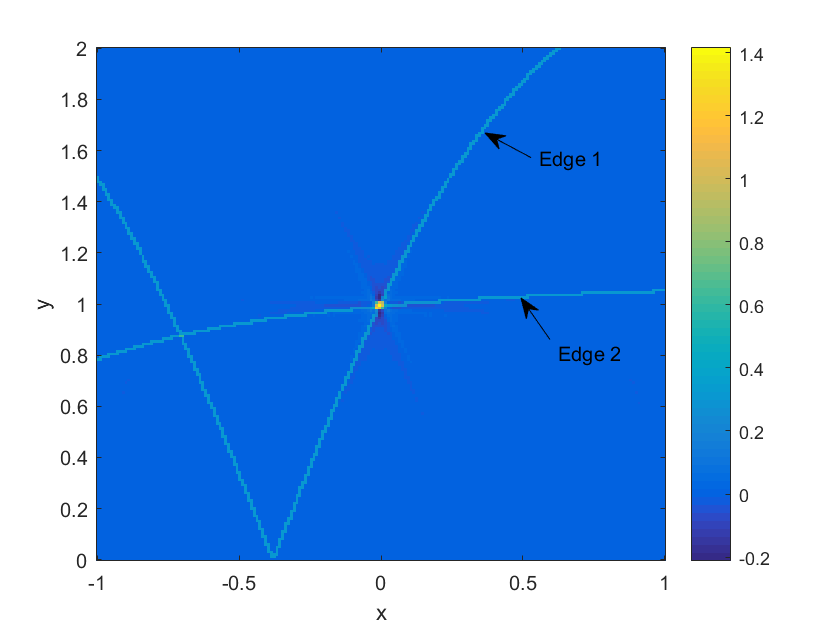}
\subcaption{Landweber.}\label{FC1:C}
\end{subfigure}
\begin{subfigure}{0.32\textwidth}
\includegraphics[width=0.9\linewidth, height=4cm]{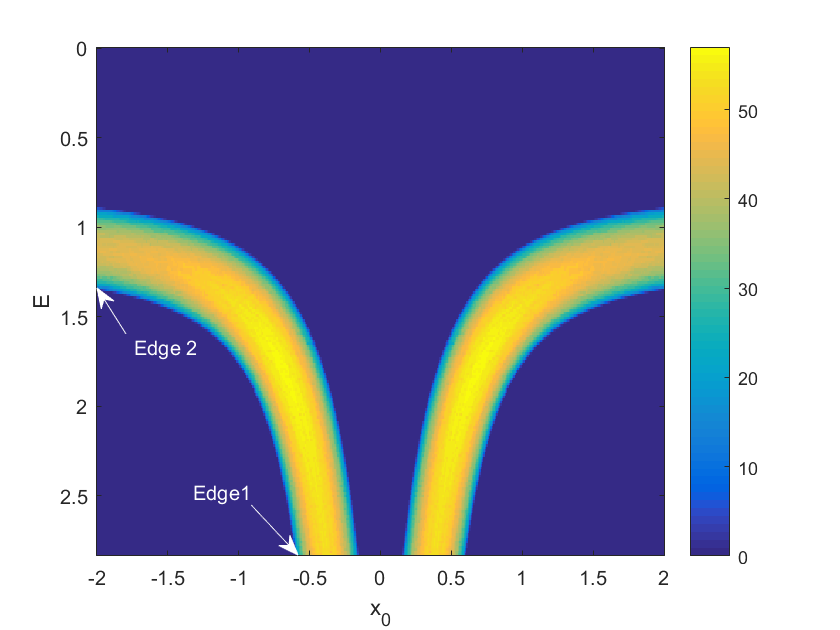}
\subcaption{$R\chi_B$ sinogram.}\label{FC1:D}
\end{subfigure}
\begin{subfigure}{0.32\textwidth}
\includegraphics[width=0.9\linewidth, height=4cm]{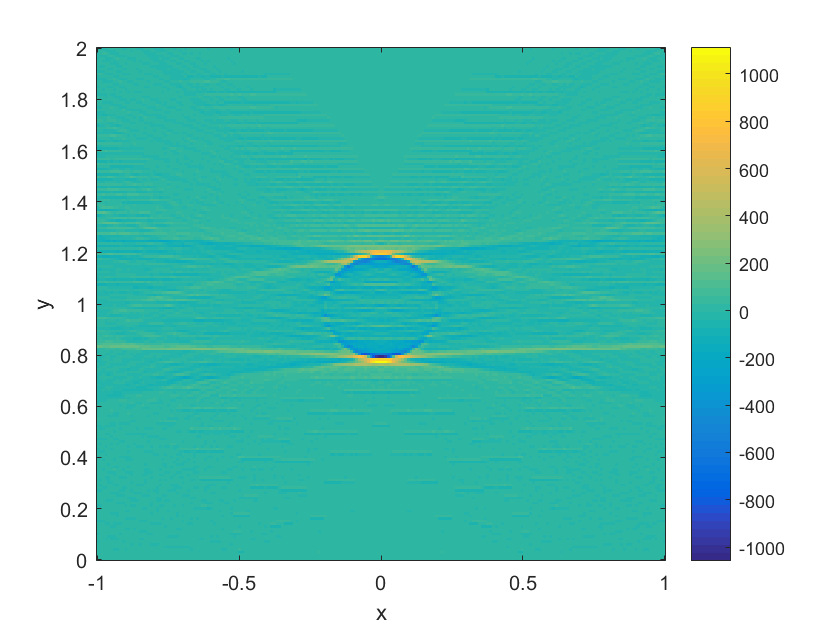}
\subcaption{$R^*\frac{\dd^2}{\dd E^2}R\chi_B$.}\label{FC1:E}
\end{subfigure}
\begin{subfigure}{0.32\textwidth}
\includegraphics[width=0.9\linewidth, height=4cm]{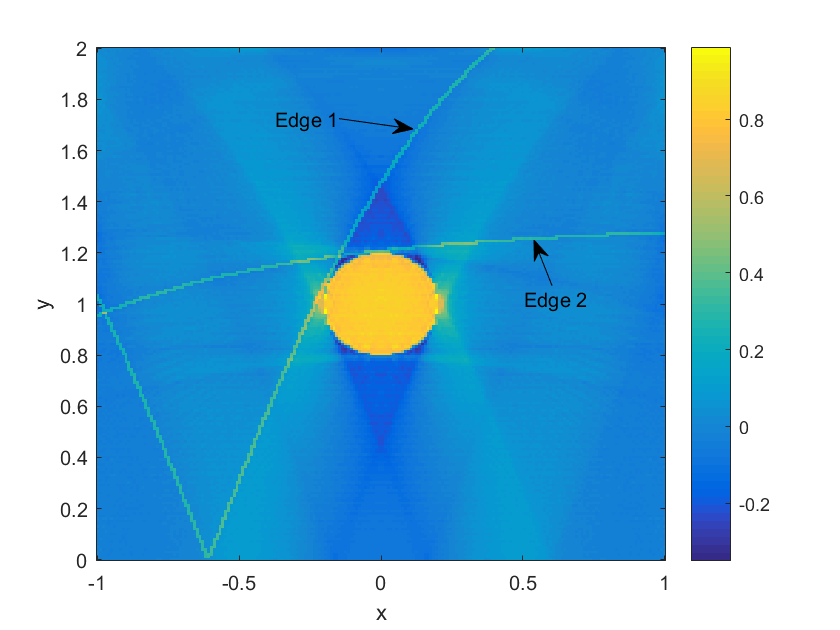}
\subcaption{Landweber.}\label{FC1:F}
\end{subfigure}
\caption{Bragg transform reconstructions. Two edges of the data set are highlighted in each case, which correspond to artefacts in the reconstruction. These are highlighted in the Landweber image.} \label{FC1}
\end{figure}
$$g'_B(r)=\frac{q''(r)}{q(r)}-\frac{q'(r)^2}{q^2(r)}=-\frac{3r^2+1}{r^2(r^2+1)^{2}}<0.$$
Thus the Bolker assumption holds for $R_j$ and $R$ when $q=q_B$ in
BST. See figure \ref{FC1} for reconstructions of $f=\delta$ and
$f=\chi_{B}$ from $Rf$, when $q=q_B$. The scanning region is
$[-1,1]\times [0,2]$ and $Rf$ is simulated for $E\in (0,2.83)$ and
$x_0\in[-2,2]$, as in example \ref{ex1}. We see artefacts appearing
along Bragg curves at the boundary of the dataset, due to the cutoff
in sinogram space. 
% The points on the sinograms in figures \ref{FC1:A} and \ref{FC1:D}
% labeled ``Edge 1" and ``Edge 2" are in the support of the data on the
% boundary of the data set, and the sharp cutoff creates artifacts as
% highlighted along Bragg curves in figures \ref{FC1:C} and \ref{FC1:F}.
The points on the sinograms in figures \ref{FC1:A} and \ref{FC1:D}
labeled ``Edge 1" and ``Edge 2" are in the support of the data and on
the boundary of the data set; the sharp cutoff at the boundary creates
artifacts as highlighted along Bragg curves (which are labeled ``Edge
1'' and ``Edge 2'' in figures \ref{FC1:C} and \ref{FC1:F}).  Similar
to example \ref{ex1}, only half of the Bragg curve at Edge 2
intersects $[-1,1]\times [0,2]$, and hence the Edge 2 artefacts appear
along one-sided Bragg curves (minus the reflected curve in $x=x_0$) in
figures \ref{FC2:C} and \ref{FC2:F}.  There is a horizontal and
vertical blurring due to limited data in the Landweber reconstruction.
This observation is in line with the theory of section
\ref{sect:Sobolev} and Theorem \ref{thm:lambda}. We noticed a similar
effect in reconstructions from broken-ray curves, when $q=q_C$ in
example \ref{ex1}. In this case the vertical blurring is less
pronounced. The vertical singularities appear sharper in the FBP
reconstructions also. This is because of the flatter gradients of the
Bragg curves as $r\to\infty$, compared to straight lines, which allow
the Bragg curves to better detect vertical singularities. That is the
Bragg curves are such that 
$$0=\lim_{r\to\infty}q'_B(r)<\min\paren{\{q'_C(r): r\in[0,\infty)\}}=1.$$
See figure \ref{Bedge} for a visualization. We display a shifted Bragg curve $q_B$ and Compton curve $q_C$. $q_B$ and $q_C$ intersect on the boundary of $\chi_B$ at $(0,1.2)$, where a singularity occurs in the direction $(0,1)$ (a vertical singularity). The gradients at $(0,1.2)$ are, $q'_C(0)=1.2$ and $q'_B(0)\approx 0.6$, approximately half the gradient of $q_C$ at $(0,1.2)$. The reduction in gradient allows for better detection of the singularity at $(0,1.2)$ using Bragg curves.
\begin{figure}[!h]
%\hspace*{-1cm}
\begin{subfigure}{0.45\textwidth}
\includegraphics[width=0.9\linewidth, height=6cm,keepaspectratio]{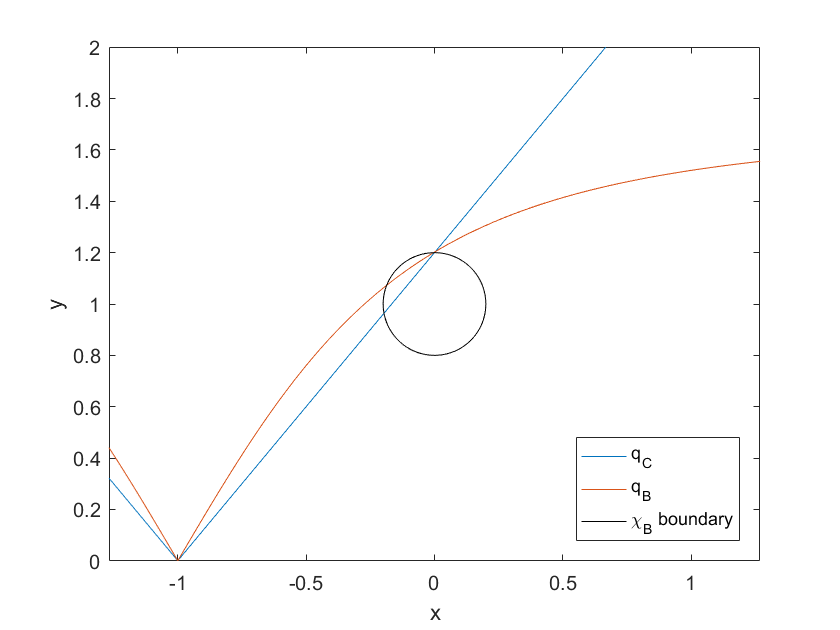}
\end{subfigure}
\caption{Illustration of Compton vs Bragg, vertical edge detection.}
\label{Bedge}
\end{figure}
\end{example}

\begin{example}[sinusoidal curves: Bolker not satisfied]
\label{ex4}
Here we give an example $q$ which satisfies \eqref{hyp1}, but fails to satisfy the Bolker assumption. We define the sinusoidal curves as
\begin{equation}
q_S(r)=(1+\epsilon)r+\sin r,
\end{equation}
where $\epsilon>0$. See figure \ref{GC4:A}.
We can check that $q_S(0)=0$ and $q'_S(r)=(1+\epsilon)+\cos r>\epsilon>0$, and hence \eqref{hyp1} is satisfied. As \eqref{hyp1} holds, it follows that $R$ is injective by \cite[Theorem 5.2]{Web}, and hence there are no artefacts due to null space. We have
$$g_S(r)=\frac{(1+\epsilon)+\cos r}{(1+\epsilon)r+\sin r},$$
which is non-injective. See figure \ref{GC4:B}.
Further $
g'_S(r)=\frac{g_1(r)}{q^2_S(r)}$,
where
$$g_1(r)=(1+\epsilon)\paren{r\sin r+2\cos r}+(1+\epsilon)^2+1,$$ and hence $g'_S(r)=0\iff g_1(r)=0$, for $r\in(0,\infty)$. $g_1$ is zero for infinitely many $r\in\mathbb{R}$, for any $\epsilon$ chosen. See figure \ref{GC4:C}.
\begin{figure}[!h]
%\centering
\begin{subfigure}{0.32\textwidth}
\includegraphics[width=0.9\linewidth, height=4cm]{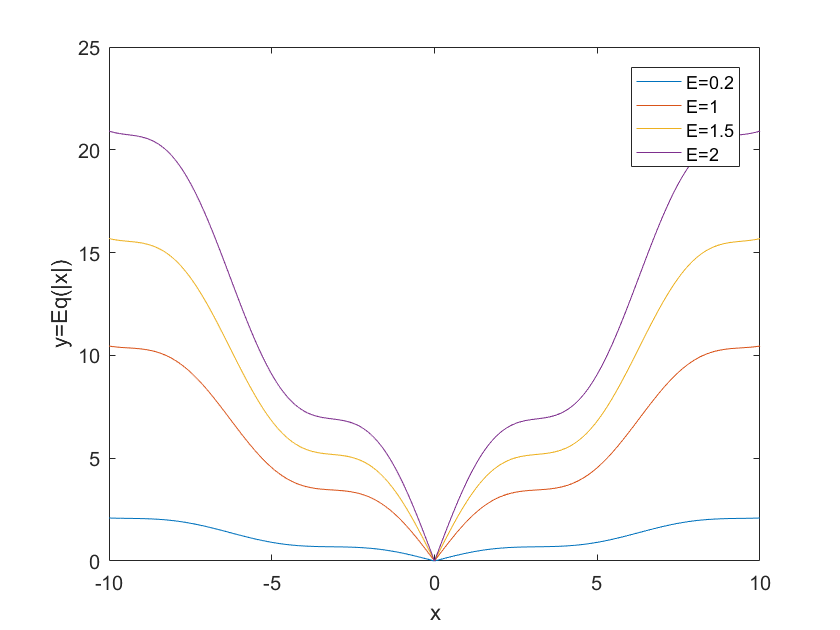}
\subcaption{$Eq_S$ for $E\in\{0.2,1,1.5,2\}$.}\label{GC4:A}
\end{subfigure}
\begin{subfigure}{0.32\textwidth}
\includegraphics[width=0.9\linewidth, height=4cm]{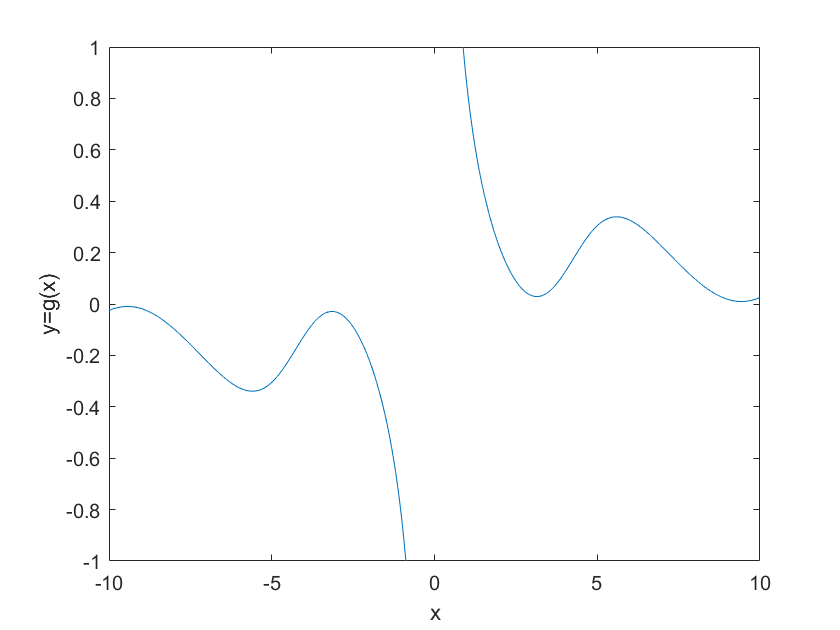}
\subcaption{$g_S$.}\label{GC4:B}
\end{subfigure}
\begin{subfigure}{0.32\textwidth}
\includegraphics[width=0.9\linewidth, height=4cm]{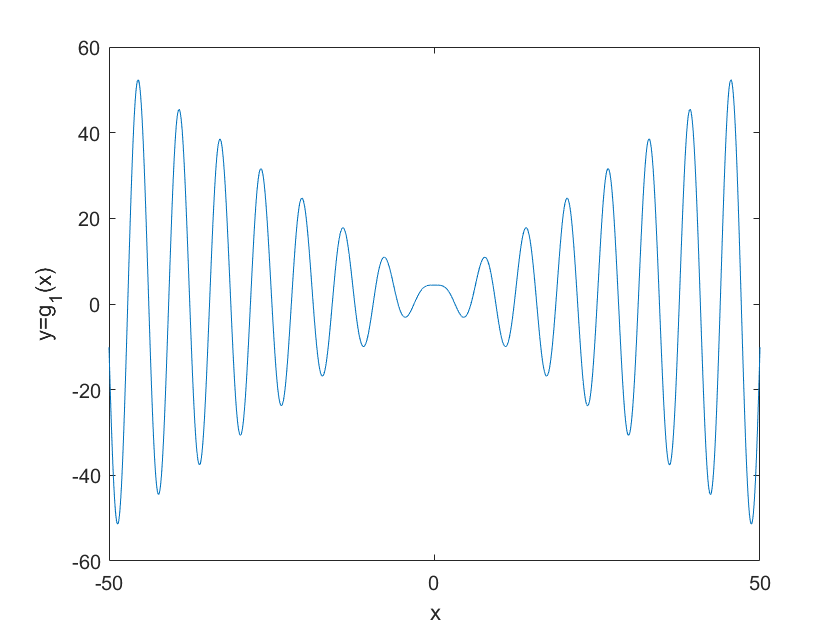}
\subcaption{$g_1$.}\label{GC4:C}
\end{subfigure}
\caption{The curves $q_S$, $g_S$ and $g_1$, for $\epsilon=0.1$. Note that the $x$ and $y$ axis limits vary across sub-figures (A), (B) and (C).} \label{GC4}
\end{figure}Hence for the sinusoidal curves the Bolker assumption is
not satisfied, and we can expect to see artefacts due to $\Lambda$
(see Remark \ref{rem3.3}) in the reconstruction. We note that $R$ and
the $R_j$ are injective by \cite[Theorem 5.2]{Web} and hence we expect
no additional artefacts due to null space. The scanning region used in
this example is $[-10,10]\times [0,20]$, and $Rf$ is simulated for
$E\in (0,3.77)$ and $x_0\in[-20,20]$. We scale up by a factor of $10$
in this case to allow for multiple oscillations of the sinusoidal
curves within the scanning region. See figure \ref{GC4:A}. On
$[-1,1]\times [0,2]$ (the scanning region used in examples \ref{ex1}
and \ref{ex2}), the $q_S$ curves are appear approximately as
broken-rays (V-lines) in the simulations, since $\sin r\approx r$ for
$r$ close to zero. Hence we scale the scanning region size by 10 here
to better highlight the discrepancies between broken-ray and
sinusoidal transform reconstruction. In these dimensions the $x_0$
range used is the same (relatively speaking) as in examples \ref{ex1}
and \ref{ex2}. The energy range is chosen so that the integration
curves have a wide variety of gradients, and sufficiently cover
$[-10,10]\times [0,20]$.
\begin{figure}[!h]
\centering
\begin{subfigure}{0.35\textwidth}
\includegraphics[width=1\linewidth, height=4.8cm]{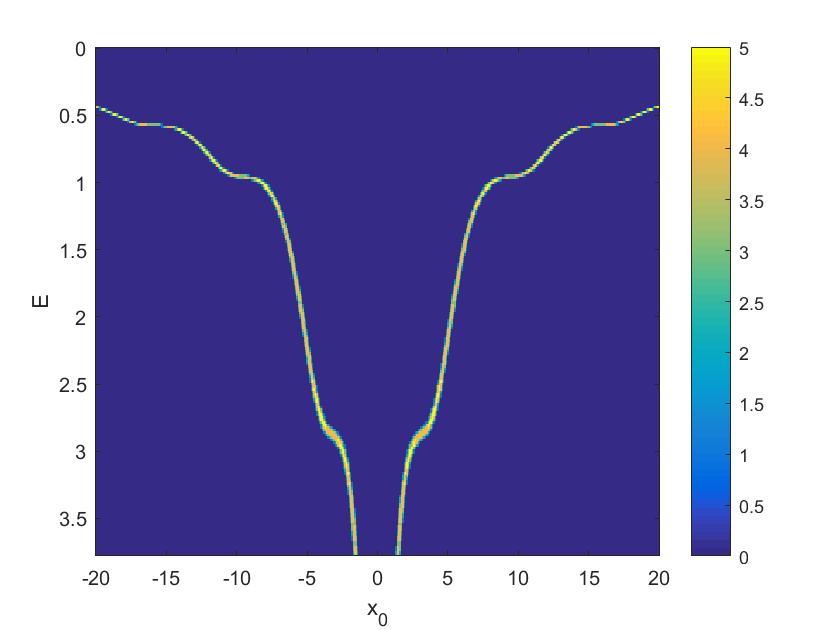} 
\subcaption{$R\delta$  sinogram.}
\end{subfigure}
\begin{subfigure}{0.35\textwidth}
\includegraphics[width=1\linewidth, height=4.8cm]{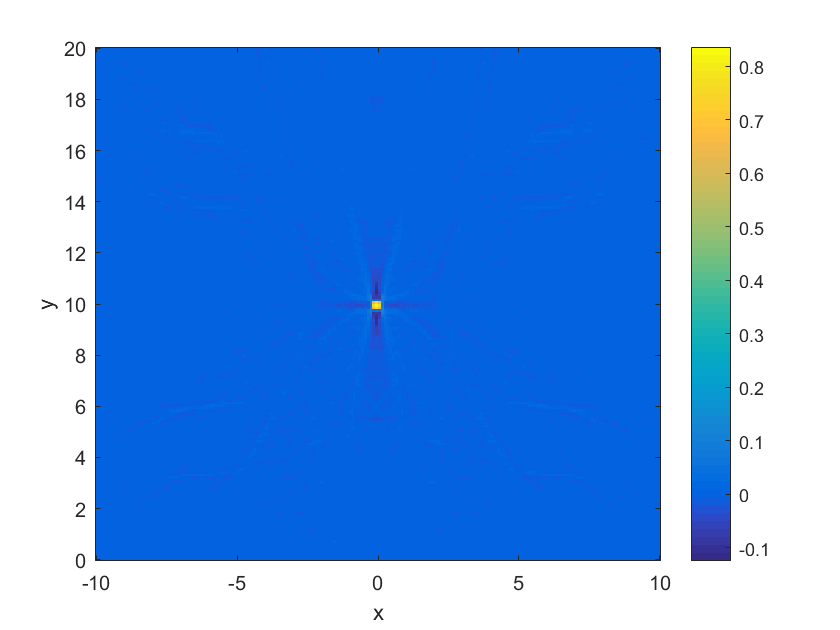}
\subcaption{Landweber.}\label{FC4:C}
\end{subfigure}
\begin{subfigure}{0.35\textwidth}
\includegraphics[width=1\linewidth, height=4.8cm]{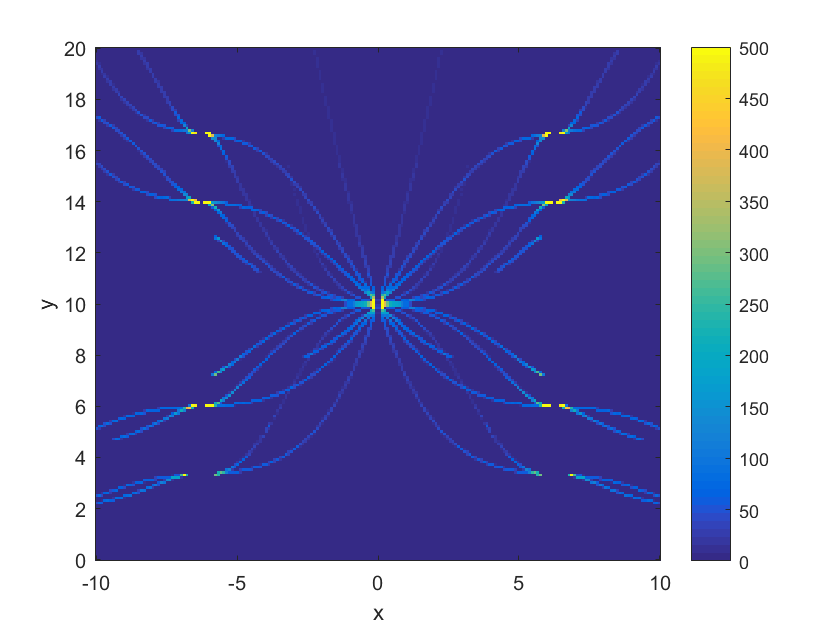} 
\subcaption{$\Lambda$  artefacts.}\label{FC4:A}
\end{subfigure}
\begin{subfigure}{0.35\textwidth}
\includegraphics[width=1\linewidth, height=4.8cm]{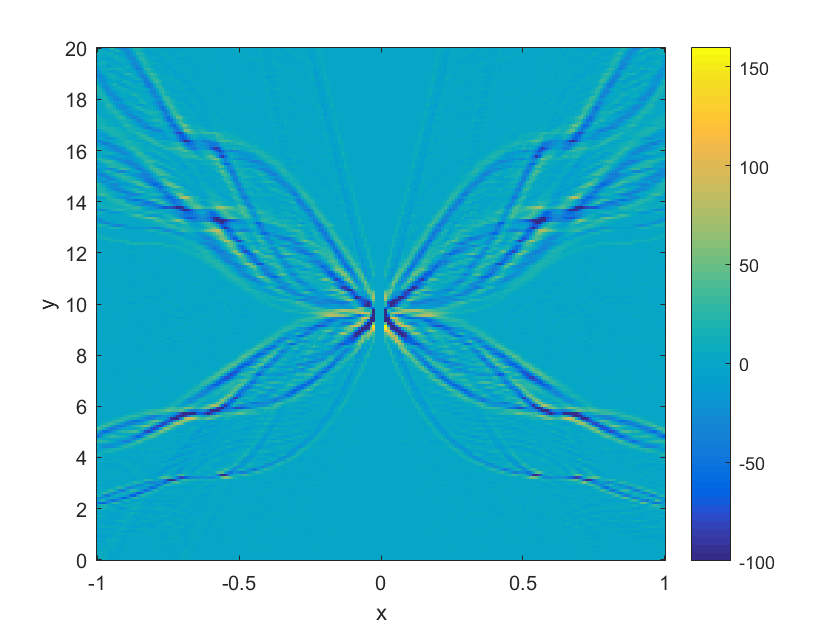}
\subcaption{$R^*\frac{\dd^2}{\dd E^2}R\delta$.}\label{FC4:B}
\end{subfigure}
\caption{Sinusoid transform reconstructions of $f=\delta$.}
\label{figPhan} \label{FC5}
\end{figure}

See figure \ref{FC5}
for reconstructions of $f=\delta$, centered at $(0,10)$, from $Rf$,
when $q=q_S$. As predicted, we see artefacts appearing in the reconstructions on the sinusoidal curves which intersect $f$ normal to a singularity (equivalently, any curve which intersects $f=\delta$). As described in Remark \ref{rem3.3}, we use $g_S$ to map the singularities of $\delta$ (at $(0,10)$, in all directions) to artefacts along sinusoidal curves. The artefacts predicted by $g_S$ and our theory are shown in figure \ref{FC4:A}. The same artefacts are observed in the FBP reconstruction in figure \ref{FC4:B}, and align exactly with our predictions. Note that we have removed the reconstructed delta function from figure \ref{FC4:B} (i.e. we set the central three image columns to zero) and truncated the color bars, to better show the artefacts. The artefacts are also observed faintly in the Landweber reconstruction in figure \ref{FC4:C}. This is in line with the theory of \cite{webber2020joint}, where the microlocal artefacts are shown to be highlighted in FBP reconstructions. 
\begin{figure}[!h]
\begin{subfigure}{0.32\textwidth}
\includegraphics[width=0.9\linewidth, height=4cm]{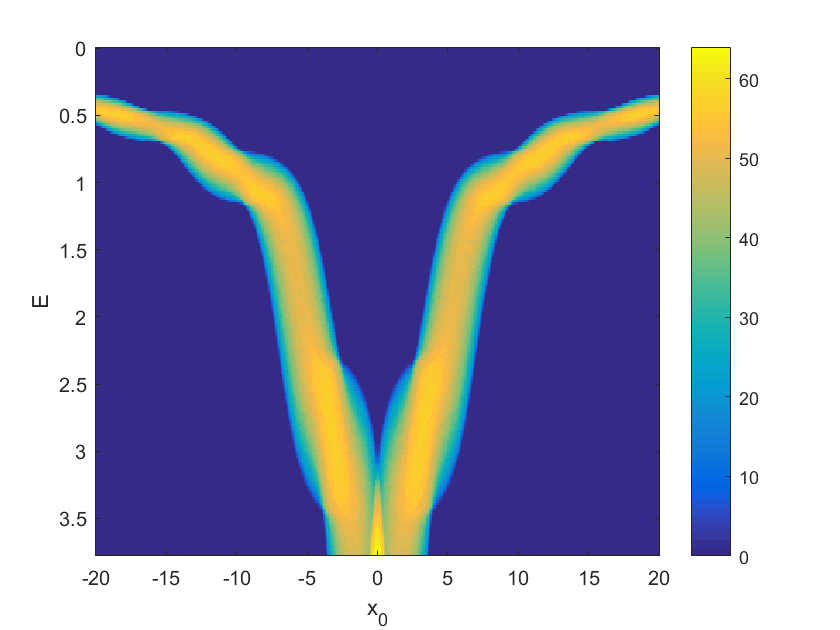}
\subcaption{$R\chi_B$ sinogram.}\label{FC4:D}
\end{subfigure}
\begin{subfigure}{0.32\textwidth}
\includegraphics[width=0.9\linewidth, height=4cm]{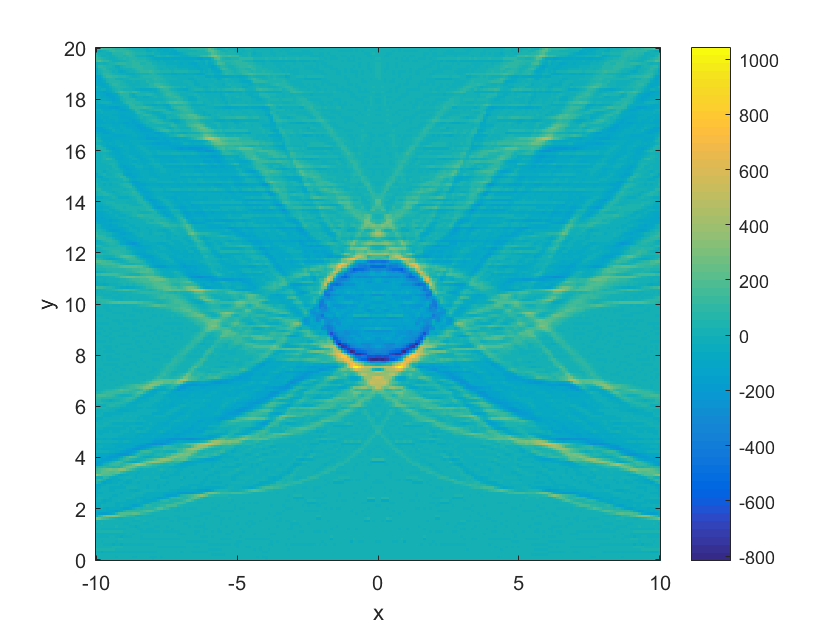}
\subcaption{$R^*\frac{\dd^2}{\dd E^2}R\chi_B$.}\label{FC4:E}
\end{subfigure}
\begin{subfigure}{0.32\textwidth}
\includegraphics[width=0.9\linewidth, height=4cm]{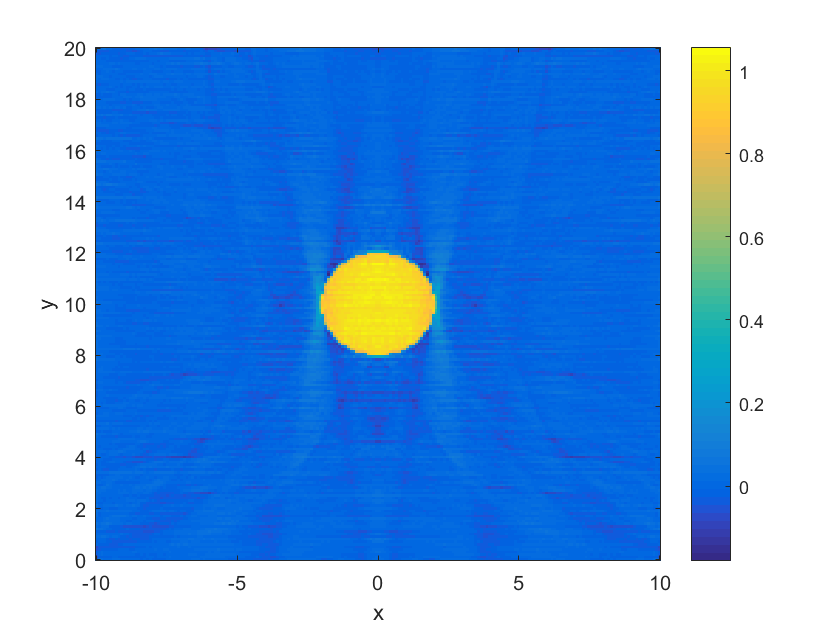}
\subcaption{Landweber.}\label{FC4:F}
\end{subfigure}
\caption{Sinusoid transform reconstructions of $f=\chi_{B}$.} \label{FC4}
\end{figure}

Let $B=\{(x,y)\in\mathbb{R}^2 : x^2+(y-10)^2<4\}$. See figure
\ref{FC4} for reconstructions of $f=\chi_{B}$ from $Rf$. Similar to
the $f=\delta$ case, we see artefacts appearing on the sinusoidal
curves which are tangent to the boundary of $\chi_B$ (i.e. the set of
points where $f$ has singularities). See figure \ref{FC4:E}. Further
$q_S$ has a larger range of gradients, when compared to $q_C$ and
$q_B$. That is
$\mu\paren{q'_S([0,\infty))}>\mu\paren{q'_C([0,\infty))},\mu\paren{q'_B([0,\infty))}$,
where $\mu$ is Lebesgue measure. Hence the Radon data can resolve the
image singularities in more directions with sinusoidal curves, when
compared to BST and CST curves. See figure \ref{Bedge} and the
arguments towards the end of example \ref{ex2}. Due to the increased
range of gradients, the horizontal and vertical blurring effects
observed in examples \ref{ex1} and \ref{ex2} are less prominent here.
This is evidenced by figure \ref{FC4:F}.
\end{example}